\documentclass[]{amsart}

\usepackage[utf8]{inputenc}
\usepackage[mathscr]{eucal}
\usepackage{stmaryrd}
\usepackage{amsmath,amsthm,amsfonts,amssymb,amscd}
\usepackage[usenames,dvipsnames]{xcolor}
\usepackage{xspace}
\usepackage{lastpage}
\usepackage{fancyhdr}
\usepackage{graphicx}
\usepackage{listings}
\usepackage{hyperref}
\usepackage{enumitem}
\usepackage{relsize}
\usepackage{tikz-cd}
\usepackage{mdwlist}
\usepackage{multicol}
\usepackage{stmaryrd}
\usepackage{float}
\usepackage{adjustbox}
\usepackage{tikz}
\usepackage{mathbbol}
\usetikzlibrary{positioning}

\usepackage[noabbrev,nameinlink]{cleveref}
\Crefformat{section}{#2\S#1#3}

\tikzset{sgplattice/.style={inner sep=1pt,norm/.style={red!50!blue},char/.style={blue!50!black},
        lin/.style={black!50}},cnj/.style={black!50,yshift=-2.5pt,left=-1pt of #1,scale=0.5,fill=white}}

\usetikzlibrary{matrix, calc, arrows}

\hypersetup{
    colorlinks=true, 
    linktoc=all,     
    linkcolor=Brown,citecolor=Brown,urlcolor=MidnightBlue,  
}

\setlist[enumerate]{label=(\alph*)}

\usetikzlibrary{graphs,decorations.pathmorphing,decorations.markings}
\tikzcdset{scale cd/.style={every label/.append style={scale=#1},
        cells={nodes={scale=#1}}}}

\newcommand{\Z}{\mathbb{Z}}
\newcommand{\R}{\mathbb{R}}

\newcommand{\N}{\mathbb{N}}
\newcommand{\F}{\mathbb{F}}

\newcommand{\Hom}{\operatorname{Hom}}
\newcommand{\Aut}{\operatorname{Aut}}

\newcommand{\tr}{\operatorname{tr}}

\newcommand{\inv}{^{-1}}

\newcommand{\End}{\operatorname{End}}
\newcommand{\op}{^{\operatorname{op}}}
\newcommand{\Res}{\operatorname{Res}}
\newcommand{\Ind}{\operatorname{Ind}}
\newcommand{\Inf}{\operatorname{Inf}}

\newcommand{\id}{\operatorname{id}}

\newcommand{\im}{\operatorname{im}}

\newcommand{\Syl}{\operatorname{Syl}}

\newcommand{\calB}{\mathcal{B}}

\newcommand{\calF}{\mathcal{F}}

\newcommand{\calL}{\mathcal{L}}
\newcommand{\calM}{\mathcal{M}}
\newcommand{\calN}{\mathcal{N}}
\newcommand{\calO}{\mathcal{O}}

\newcommand{\calS}{\mathcal{S}}

\newcommand{\calY}{\mathcal{Y}}

\newcommand{\comp}{\on{comp}}

\newcommand{\catC}{\mathscr{C}}
\newcommand{\catD}{\mathscr{D}}

\newcommand{\catI}{\mathscr{I}}

\newcommand{\catK}{\mathscr{K}}
\newcommand{\catL}{\mathscr{L}}
\newcommand{\catM}{\mathscr{M}}

\newcommand{\catP}{\mathscr{P}}
\newcommand{\catQ}{\mathscr{Q}}

\newcommand{\frakp}{\mathfrak{p}}
\newcommand{\frakq}{\mathfrak{q}}

\newcommand{\catmod}{\on{mod}}

\newcommand{\Spc}{\operatorname{Spc}}
\newcommand{\Spec}{\operatorname{Spec}}
\newcommand{\Pic}{\operatorname{Pic}}
\newcommand{\Weyl}[2]{{#1}/\!\!/{#2}}
\newcommand{\open}{\operatorname{open}}
\newcommand{\cone}{\operatorname{cone}}
\newcommand{\catproj}{\mathsf{proj}}
\newcommand{\catperm}{\mathsf{perm}}
\newcommand{\on}[1]{\operatorname{#1}}
\newcommand{\supp}{\on{supp}}
\newcommand{\sbull}{{\scriptscriptstyle\bullet}}
\newcommand{\tH}{\on{H}^{\sbull\sbull}}

\newcommand{\Sub}{\on{Sub}}

\makeatletter
\newcommand{\xhookdoubleheadrightarrow}[2][]{%
    \lhook\joinrel
    \ext@arrow 0359\rightarrowfill@ {#1}{#2}%
    \mathrel{\mspace{-15mu}}\rightarrow
}
\makeatother

\newtheorem{theorem}{Theorem}[section]
\newtheorem{lemma}[theorem]{Lemma}
\newtheorem{prop}[theorem]{Proposition}
\newtheorem{corollary}[theorem]{Corollary}

\newtheorem*{theorem*}{Theorem}

\newtheorem*{conjecture*}{Conjecture}

\theoremstyle{remark}
\newtheorem{remark}[theorem]{Remark}
\newtheorem{example}[theorem]{Example}

\newtheorem*{question*}{Question}

\theoremstyle{definition}
\newtheorem{definition}[theorem]{Definition}
\newtheorem{notation}[theorem]{Notation}
\newtheorem{construction}[theorem]{Construction}
\newtheorem*{definition*}{Definition}

\newtheorem{innercustomthm}{Theorem}
\newenvironment{customthm}[1]
{\renewcommand\theinnercustomthm{#1}\innercustomthm}
{\endinnercustomthm}

\AddToHook{env/corollary/begin}{\crefalias{theorem}{corollary}}
\AddToHook{env/prop/begin}{\crefalias{theorem}{proposition}}
\AddToHook{env/lemma/begin}{\crefalias{theorem}{lemma}}
\AddToHook{env/definition/begin}{\crefalias{theorem}{definition}}
\AddToHook{env/example/begin}{\crefalias{theorem}{example}}
\AddToHook{env/remark/begin}{\crefalias{theorem}{remark}}
\AddToHook{env/observation/begin}{\crefalias{theorem}{observation}}
\AddToHook{env/construction/begin}{\crefalias{theorem}{construction}}
\AddToHook{env/recollection/begin}{\crefalias{theorem}{recollection}}
\AddToHook{env/notation/begin}{\crefalias{theorem}{notation}}

\begin{document}
    \title{Permutation twisted cohomology, remixed}
    \author{Sam K. Miller}
    \date{September 23, 2026}
    \address{Department of Mathematics, University of Georgia, Athens GA, 30602} 
    \email{sam.miller@uga.edu} 
    \urladdr{https://www.samkmiller.com/}
    \begin{abstract}
        For each endotrivial complex for a $p$-group arising from Bredon homology of a representation sphere, we construct $p$-local quasi-isomorphisms, called forerunners. These enable us to extend Balmer--Gallauer's results in \cite[Part II]{BG25} concerning the tensor-triangular geometry of permutation modules for elementary abelian $p$-groups to all $p$-groups. We construct an open cover of the Balmer spectrum under which all endotrivials are tt-line bundles, that is, every endotrivial is locally isomorphic to a shift of the tensor unit. We define a remixed permutation twisted cohomology ring for which the canonical comparison map from the Balmer spectrum to the homogeneous spectrum of the twisted cohomology ring is injective. If the twisted cohomology ring is Noetherian, the comparison map is an open immersion, and the open cover endows the Balmer spectrum with Dirac scheme structure. We prove Noetherianity holds for Dedekind groups and the dihedral group of order 8, and conjecture it holds for all $p$-groups.
    \end{abstract}

    \maketitle
    For the duration of this paper, $k$ denotes a field of prime characteristic $p$ and $G$ denotes a finite group.

    \section*{Introduction}
    It is oft claimed that \emph{permutation $kG$-modules}, modules which admit a $G$-stable $k$-basis, and their direct summands, \emph{$p$-permutation $kG$-modules}, are at minimum the second easiest class of representations of $G$ over $k$ to define. In a wildly (c.f. \cite{BD77}) vast representation-theoretic cosmos of $kG$-modules, our $p$-permutation modules appear tiny and insignificant, drowning in the endless sea of symmetries. However, our protagonist, the bounded homotopy category of $p$-permutation $kG$-modules $\catK(G) := \on{K}_b(\mathsf{perm}(kG)^\natural)$, exhibits significant control over said cosmos. For a classical example, if $G$ has an abelian Sylow $p$-subgroup $S$, Brou\'{e}'s abelian defect group conjecture \cite{Bro90, R96} predicts a particular triangulated equivalence $\on{K}_b(\mathsf{perm}(kGb_0)^\natural) \cong \on{K}_b(\mathsf{perm}(kN_G(S)c_0)^\natural)$ where $b_0$ and $c_0$ denote the principal blocks of $kG$ and $kN_G(S)$ respectively; such an equivalence induces a corresponding derived equivalence. More recently, Balmer--Gallauer demonstrated that every $kG$-module has a finite resolution by $p$-permutation modules \cite{BG23}. Consequently, $\catK(G)$ admits the bounded derived category $\on{D}_b(\mathsf{mod}(kG))$ as a Verdier quotient, strengthening an unpublished result of Rouquier \cite{Rou06}. Although these $p$-permutation modules are few, they are mighty; together, they recover the entire representation theory of $kG$!

    In a landmark series of papers \cite{BG22, BG23, BG23b, BG25}, Balmer--Gallauer elucidated the tensor-triangular geometry of $\catK(G)$. In particular, in \cite{BG25}, the authors deduced its Balmer spectrum $\Spc(\catK(G))$, the universal support of the category.\footnote{This classification has ramifications beyond representation theory: to the geometrically or topologically inclined reader, this classification can be interpreted in terms of cohomological Mackey functors, Artin motives, or the homotopy category of modules over the constant Green functor, see \cite{BG23b, F25}} As a set, the Balmer spectrum is built from $p$-local cohomology: \[\Spc(\catK(G)) = \bigsqcup_{H \in \Sub_p(G)/G} \Spc(\on{D}_b(k[\Weyl{G}{H}])) \cong \bigsqcup_{H \in \Sub_p(G)/G} \Spec^h(\on{H}^\sbull(\Weyl{G}{H};k)).\] This stratification is obtained from certain tt-functors classically known as \emph{Brauer quotients} for modules, and in \cite{BG25} are termed \emph{modular fixed points}, functors $\Psi^H\colon \catK(G) \to \catK(\Weyl{G}{H})$, where $\Weyl{G}{H} := N_G(H)/H$ is the \emph{Weyl group} of $G$ at $H$. These functors mimic the fixed-point operation for $G$-sets: given any $G$-set $X$ and $p$-subgroup $H \leq G$, one has a natural isomorphism $\Psi^H(kX) \cong k[X^H]$. We'll discuss the nitty-gritty of the classification, but first take a narrative detour to the domain of endotrivials.

    \subsection*{Endotrivial complexes and forerunners} Separately, in a sequence of papers \cite{M24a, M24b, M24c}, we deduced the Picard group $\Pic(\catK(G))$, i.e., the group of \emph{endotrivial complexes} for all finite groups. Endotrivials are by definition the invertible objects of $\catK(G)$, the chain complexes $C \in \catK(G)$ satisfying $C^* \otimes C \simeq k[0]$. This classification was also partially determined in Bachmann's dissertation \cite{Bac16} in the context of Artin motives. Independently, a similar computation was performed by Yal\c{c}in \cite{Y16} for $p$-groups in the context of $G$-Moore spaces, which was partially extended to finite groups in joint work with Gelvin \cite{GY21}.

    There are multiple alternative ways to think about endotriviality that lead to a classification. In \cite{M24b, M24c}, we deduced $\Pic(\catK(G))$ by considering endotrivials as a special class of more generally `relatively' endotrivial complexes akin to \emph{endopermutation modules} (see \cite[Chapter 12]{Bou10}); this parallel to endopermutation modules is developed further in \cite{M26a}. Another characterization considered by Bachmann and Gelvin--Yal\c{c}in is topological in nature: endotrivials arise from Bredon homology of (virtual) \emph{representation spheres}, i.e., one-point compactifications of real representations. In other words, one has a homomorphism $\on{RO}(G) \to \Pic(\catK(G))$. In fact, for $p$-groups (but not in general), the homomorphism is surjective, and its kernel is characterized by the Adams operation or the $j$-homomorphism. In particular, for $G$ a $p$-group, $\Pic(\catK(G))$ has a canonical $\Z$-basis associated to the irreducible real representations of $G$.

    In this paper, we adopt this topological point of view for $p$-groups, with a particular emphasis on \emph{effective} endotrivials, the endotrivials arising from Bredon homology of (genuine) representation spheres. Such endotrivials have additional structure: they are chain complexes of free modules over the \emph{orbit category} $\Gamma_G$. This structure allows us to construct certain $p$-local quasi-isomorphisms which we term \emph{forerunners}, as they foretell quasi-isomorphisms prior to applying modular fixed points.

    \begin{customthm}{A}(\Cref{thm:construction})
        Let $G$ be a $p$-group. For every effective endotrivial complex $C$ (see \Cref{def:decreasing}) and subgroup $H \leq G$, there exists a \emph{forerunner} homomorphism $\iota^H_C: k \to C[-h_C(H)]$ such that $\Psi^H(\iota^H_C)$ is a quasi-isomorphism.
        Here, $h_C$ denotes the h-mark homomorphism of $C$, see \Cref{rmk:hmarks}.
    \end{customthm}

    Let us justify to the inquisitive reader why the existence of forerunners is noteworthy. First, the existence of such morphisms is by no means a given: it is a special property of effective endotrivials. Proving their existence requires developing the theory of chain complexes of permutation $kG$-modules arising from the orbit category $\Gamma_G$; we do not list the results here, but believe they will be of independent interest. Second, Gallauer has already utilized these morphisms in his recent initiation of the study of \emph{periodicity} for symmetric monoidal $\infty$-categories \cite{Gal25}; we remark on this work in the sequel. In particular, using forerunners, we were able to confirm in joint work with Gallauer a conjecture he had previously posed. But the real justification, the true power of these morphisms, is that they enable us to explain the geometry of permutation modules in a much more general light, as we'll now explain.

    \subsection*{A new open cover of the spectrum} In \cite{BG25}, Balmer--Gallauer split their deduction of $\Spc(\catK)$ into two parts: Part I concerns deducing $\Spc(\catK(G))$ as a set and constructing the modular fixed-point functors, and Part II deduces the topology of $\Spc(\catK(G))$. In Part II, the authors prove a colimit theorem a la Quillen stratification reducing classification to the case of $G$ elementary abelian, then investigate this case in greater depth. Using the classification of endotrivials, the existence of forerunners, and the topological origins of effective endotrivials, we extend these investigations in the context of any finite $p$-group. Our hope is that $\Spc(\catK(G))$ can be completely deduced without any gluing. The $p$-group case is the crux: given a finite group $G$ with $p$-Sylow $S$, one can determine $\Spc(\catK(G))$ from $\Spc(\catK(S))$ by restriction via \cite[Propositions 3.8, 4.7]{BG25}.

    We first review Balmer--Gallauer's methods and results. Let $E$ be an elementary abelian $p$-group. Balmer--Gallauer constructed an open cover $\{U(H)\}_{H \leq E}$ of $\Spc(\catK(E))$ indexed by subgroups of $E$. Explicitly, these opens are obtained by considering the open loci of certain morphisms from the tensor unit to certain endotrivials. This cover has some remarkable features:

    \begin{enumerate}
        \item All endotrivials for $E$ are isomorphic to a shift of the tensor unit in the localization $\catK(E)|_{U(H)}$; in other words, all endotrivials for $E$ are \emph{line bundles} \cite[Remark 13.10]{BG25};
        \item Each open $U(H)$ contains the closed point $\catM(H) \in \Spc(\catK(E))$. Hence $\{U(H)\}_{H \leq E}$ is an open cover of $\Spc(\catK(E))$; \cite[Proposition 13.11]{BG25}
        \item Under this open cover, there is an identification $\catK(E)|_{U(1)} \cong \on{D}_b(kE)$. In other words, $U(1)$ is nothing more than the \emph{cohomological open}, corresponding to the image under $\Spc$ of the localization $\catK(E) \twoheadrightarrow \on{D}_b(kE)$; \cite[Proposition 13.14]{BG25}
    \end{enumerate}

    This open cover endows $\Spc(\catK(E))$ with \emph{Dirac scheme structure} in the sense of \cite{HP23}; a Dirac scheme is to a usual scheme as a graded ring is to a non-graded one. The morphisms used to construct these opens are examples of forerunners - this observation suggests that such an open cover in fact exists in much greater generality.

    \begin{customthm}{B}
        Let $G$ be a $p$-group. We have open covers
        $\{U(H)\}_{H \in \Sub_p(G)/G}$ and $\{U(H)'\}_{H \in \Sub_p(G)/G}$ of $\Spc(\catK(G))$ with $U(H)' \subseteq U(H)$ (\Cref{con:open}) indexed by conjugacy classes of subgroups of $G$  for which the following hold.
        \begin{enumerate}
            \item Every endotrivial $C$ is a \emph{tt-line bundle}, i.e., over each open $U(H)$ or $U(H)'$ is isomorphic to a shift of the tensor unit in the localization $\catK(G)|_{U(H)}$. Explicitly, $C \cong k[h_C(H)]$ in $\catK(G)|_{U(H)}$, where $h_C$ denotes the h-mark homomorphism associated to $C$; (\Cref{cor:linebundle})
            \item Each open $U(H)'$ contains a \emph{unique} closed point $\catM(H) \in\Spc(\catK(G))$, corresponding to the unique closed point of $\Spec^h(\on{H}^\sbull(\Weyl{G}{H};k))$; (\Cref{cor:onlyclosedpoint})
            \item The open $U(1) = U(1)'$ is the \emph{cohomological open} $ \Spc(\on{D}_b(kG))\subseteq\Spc(\catK(G))$. Therefore, we have a tt-equivalence $\catK(G)|_{U(1)} \cong \on{D}_b(kG)$ (\Cref{prop:openisac}), and realize $\on{H}^\sbull(G)$ as the twist-zero part of a localization of $\tH(G)$ (\Cref{cons:localizations}).
        \end{enumerate}
    \end{customthm}

    \subsection*{(Remixed) permutation twisted cohomology and Dirac geometry} To completely deduce the topology of $\Spc(\catK(E))$, Balmer--Gallauer define the \emph{permutation twisted cohomology ring} $\tH(G)$. Explicitly, the authors consider morphisms from the tensor unit to shifts of certain endotrivials corresponding to subgroups $N \triangleleft G$ of index $p$. The collection of these morphisms defines the multigraded ring $\tH(G)$, and extends the usual notion of the cohomology ring of a triangulated category. The ring comes with a canonical \textit{comparison map}, \[\comp_G\colon  \Spc(\catK(G)) \to \Spec^h(\tH(G)),\] extending the usual comparison map of \cite{Ba10}. A twisted variant of the cohomology ring is necessary: since $\catK(G)$ is a homotopy category of chain complexes, the usual cohomology ring $\End^\sbull_{\catK(G)}(k)$ is a tad drab on its own. The miracle of this construction is:

    \begin{theorem*}{\cite[Theorem 10.5]{BG25}}
        Let $E$ be an elementary abelian $p$-group. The comparison map $\comp_E\colon  \Spc(\catK(E)) \to \Spec^h(\tH(E))$ identifies $\Spc(\catK(E))$ with an open subspace of the homogeneous spectrum of $\tH(E)$. In particular, $\Spc(\catK(E))$ is a Dirac scheme.
    \end{theorem*}

    Hausmann--Schwede have further considered the twisted cohomology ring for elementary abelian $2$-groups in the context of Bredon homology \cite{HS25}. The authors determined a minimal generating set and proved that for $p = 2$, $\tH(E)$ is nilpotent-free.
    Interestingly, for elementary abelian $p$-groups, the endotrivials involved in the construction of twisted cohomology are precisely the effective endotrivials we've considered, motivating our following `remixed' construction.

    \begin{definition*}[\Cref{def:twistedcohomology}]
        Let $G$ be a $p$-group. Denote the subset of the canonical $\Z$-basis of $\Pic(\catK(G))$ arising from \textit{nontrivial} irreducible real representations of $G$ by $\calB(G)$. Let $\N^{\calB(G)} = \{q\colon  \calB(G) \to \N\}$ be the \textit{monoid of twists}, i.e., tuples of non-negative integers indexed by the set $\calB(G)$. Then the $(\Z\times \N^{\calB(G)})$-graded ring \[\tH(G) = \tH(G;k) := \bigoplus_{s\in \Z}\bigoplus_{q \in \N^{\calB(G)}} \Hom_{\catK(G)} \left(k , \bigotimes_{C \in \calB(G)} C^{\otimes q(C)}[s]\right)\] is the \emph{(remixed) twisted cohomology ring}.
    \end{definition*}

    This ring may also be considered the `effective cone' of $\on{RO}(G)$-graded Bredon homology, see \cite{HS25} and the remarks within. We twist only by shifts of effective endotrivials rather than all endotrivials, as otherwise, the graded ring becomes unmanageably large (c.f. \cite[Remark 12.24]{BG25}). Analogously, the group cohomology ring $\on{H}^\sbull(G;k)$ is Noetherian and recovers the Balmer spectrum of the derived category $\on{D}_b(kG)$, whereas the Tate cohomology ring $\hat{\on{H}}{}^\sbull(G;k)$ is, for sufficiently interesting $G$, non-Noetherian and contains only nilpotents in negative degrees \cite{BC92}.

    \begin{customthm}{C}
        Let $G$ be a $p$-group. The comparison map \[\comp_G\colon \Spc(\catK(G)) \to \Spec^h(\tH(G)),\, \catP\mapsto \left\langle f \in \tH(G) \text{ homog.} \mid \cone(f) \not\in \catP\right\rangle\] is injective (\Cref{thm:injectivity}) with dense image (\Cref{cor:denseimage}). If $\tH(G)$ is Noetherian, the comparison map is an open immersion. (\Cref{cor:homeo})
    \end{customthm}

    To show this, we require certain `faithful' endotrivials which Balmer--Gallauer's construction may not have for larger groups, see \Cref{rmk:needallinvs}. In the proof however, we mainly rely on the existence of certain large endotrivials, rather than all effective endotrivials. This suggests a possible `reduced' version of the twisted cohomology ring, which we will investigate in future work.

    Since every endotrivial is a tt-line bundle with respect to our open cover $\{U(H)\}$ of Theorem B, for every open $U(H)$, the twisted cohomology ring locally identifies with the usual cohomology ring in the localization $\catK(G)|_{U(H)}$. This allows us to use the classical comparison map of \cite{Ba10} locally in each open. Morally, twisted cohomology globally describes the local data in each open. If the comparison map restricts on each open to a local homeomorphism, then we can endow $\Spc(\catK(G))$ with Dirac scheme structure, as was the case for elementary abelian $p$-groups.

    \begin{customthm}{D}
        Let $G$ be a $p$-group. If $\tH(G)$ is Noetherian, the comparison map is an open immersion and for each subgroup $H \leq G$, restricts to a local homeomorphism \[\comp_{\calL_G(H)}\colon U(H) \xrightarrow{\sim} \Spec^h(\End_{\catL_G(H)}^\sbull(k)).\] In this case $(\Spc(\catK(G)), \calO_G^\sbull)$ is a Dirac scheme. (\Cref{cor:homeo}, \Cref{cor:dirac}) Noetherianity holds for the following groups:
        \begin{enumerate}
            \item The dihedral group of 8 elements $D_8$; (\Cref{thm:d8noeth})
            \item \emph{Dedekind} $p$-groups, that is, $p$-groups with every subgroup normal. (\Cref{thm:resontonoeth})
        \end{enumerate}
    \end{customthm}

    Finally, we briefly remark on the connection to \cite{Gal25}. Gallauer considers twisted cohomology in an abstract context, and shows that when things behave particularly well, these constructions recover the Balmer spectrum and can be used to relate periodicity data between a tt-category and an associated twisted cohomology ring. The difficulty however, is showing such things behave well, and that is what we do here. In the language of \cite{Gal25}, finite generation of $\tH(G)$ implies that the maps into the effective endotrivials, which Gallauer terms \emph{global sections}, generate the topology of the spectrum. This determines a finite \emph{ample} family of tt-line bundles, which can be taken as the gradings which arise for a set of generators of $\tH(G)$. This notion of ampleness almost recovers the notion of ampleness in algebraic geometry. If this holds, $\catK(G)$ is a \emph{divisorial} tt-category. We precisely characterize the global sections for any effective endotrivial.

    \begin{customthm}{E}(\Cref{prop:identifyingkhoms})
        Let $G$ be a $p$-group, $C$ be an indecomposable effective endotrivial complex, and $s$ an integer. Then \[\Hom_{\catK(G)}(k , C[-s]) = \Hom_{\on{Ch}(\mathsf{perm}(kG)^\natural)}(k , C[-s]) = \Hom_{\mathsf{mod}(kG)}(k, C_s).\]
    \end{customthm}

    \subsection*{Open questions} Presently, the most pertinent question is whether the comparison map is an open immersion. This would be implied by the possibly stronger condition of Noetherianity of the remixed twisted cohomology ring $\tH(G)$. We conjecture both hold.

    \begin{conjecture*}
        For all $p$-groups $G$, the following hold.
        \begin{enumerate}
            \item Theorem D holds independently of the assumption of Noetherianity of $\tH(G)$;
            \item The twisted cohomology ring $\tH(G)$ is Noetherian.
        \end{enumerate}
    \end{conjecture*}

    Balmer--Gallauer show Noetherianity holds for $G$ elementary abelian \cite[Lemma 12.12]{BG25}, and in Section 7 we extend their strategy to show it holds in a wider range of cases. Our proofs of Noetherianity rely heavily on the topological origin of effective endotrivials, but unfortunately also rely on information about the classification of endotrivials specific to these groups. We emphasize that showing Noetherianity of $\tH(G)$ is likely not the only way to show that the comparison map is an open immersion, and expect there are alternate routes.

    While our twisted cohomology ring may be sufficient to capture the subtleties of $\Spc(\catK(G))$, it is obviously cumbersome. An explicit presentation in terms of generators and relations, akin to the computations in \cite[Section 17]{BG25} are at present beyond our reach for arbitrary groups, but should be more feasible for the groups for which we've shown Noetherianity of $\tH(G)$, as one now has an explicit list of generators. Additionally, a possible reduced version of $\tH(G)$, in which one twists by fewer endotrivials, may make such an analysis easier.

    Beyond twisted cohomology, there remain numerous questions related to endotrivials and permutation modules. Perhaps the most immediate concerns the forerunner maps: our construction is ad-hoc and depends on a choice of cellularization of a representation sphere. We expect the existence of a `canonical' construction or description. Such a construction would resolve whether certain expected properties for forerunners hold, and elucidate the open covers we construct. See \Cref{rmk:careiniota} for further discussion.

    Our understanding of endotrivials, effective or not, is still incomplete, and this seems to be the limiting factor for proving Noetherianity. Though endotrivials are classified, explicitly writing them down for larger groups (e.g. non-abelian with $p$-rank at least 3) is a significant challenge. One fairly elementary question is as follows: in \cite{M24c} it was shown that $\Pic(\catK(G))$ is a \emph{biset functor} in the sense of Bouc \cite{Bou10}, and `induction' of endotrivials is induced on the topological level by the norm map. How this presents on the level of $\catK(G)$ seems difficult to determine, as it is (perhaps surprisingly) not given by tensor induction of chain complexes.

    It is reasonable to ask if a generalization of twisted cohomology exists for all finite groups. Not all endotrivials $C$ arise from representation spheres, so our forerunners cannot always be constructed outside of the $p$-group setting. Moreover, if $G/[G,G]$ is not a $p$-group, then not all endotrivials for $kG$ are tt-line bundles. Explicitly, any shift of a $k$-dimension 1 $kG$-module $k_\omega$, with $\omega \in \Hom(G, k^\times)$, is endotrivial. Related, we do not know a description of the subgroup of $\Pic(\catK(G))$ generated by endotrivials induced from representation spheres; we expect this subgroup should be close to full-rank. We expect that one can construct an analogous open cover of $\Spc(\catK(G))$ from this subgroup, and that some of the techniques of this paper carry over without much difficulty.

    Finally, one may ask questions regarding descent with respect to the open cover $\{U(H)\}$ of $\Spc(\catK(G))$. In particular, there should be a descent spectral sequence associated to $\Pic(\catK(G))$. Although we prove every endotrivial is a tt-line bundle, it is currently unknown and an interesting question whether $\Pic(\catK(G)|_{U(H)}) \cong \Z$ for all $H \leq G$; we expect this to be the case. In correspondence, Gallauer gave a proof sketch of this for $p = 2$ and $G = V_4 = C_2 \times C_2$, but not much is known beyond this case.

    \subsection*{Organization} The paper is organized as follows. Section 1 reviews tensor-triangular geometry, $p$-permutation modules and modular fixed points. Section 2 reviews endotrivial complexes and Borel-Smith functions. Section 3 covers chain complexes of free modules over the orbit category, effective endotrivials, and proves some technical results. Section 4 constructs the forerunners $\iota^H_C$ and the open cover $\{U(H)\}$ of the Balmer spectrum $\Spc(\catK(G))$. Section 5 constructs the twisted cohomology ring $\tH(G)$ and proves compatibilities. Section 6 proves injectivity of the comparison map and corollaries. Finally, Section 7 proves Noetherianity of twisted cohomology for certain cases.

    \subsection*{Notation} Our notation mostly follows \cite{BG25}, which differs considerably from \cite{M24a, M24b, M24c}. The \textit{Weyl group} $\Weyl{G}{H}$ of $H \leq G$ is the subquotient $N_G(H)/H$ of $G$. Given two subgroups $H, K \leq G$, we write $H =_G K$ if $H, K$ are $G$-conjugate, and $H \leq_G K$ to denote $H$ is a subgroup of a $G$-conjugate of $K$. We denote the set of all subgroups of $G$ by $\Sub(G)$, the set of all $p$-subgroups of $G$ by $\Sub_p(G)$, and $\Sub_p(G)/G$ for conjugacy classes of $p$-subgroups. We write $\on{CF}(G)$ to denote the additive group of \emph{superclass functions}, i.e., functions $\Sub(G) \to \Z$ constant on conjugacy classes, and if $G$ is not a $p$-group, we write $\on{CF}(G,p)$ to denote the group of superclass functions valued on $p$-subgroups of $G$.

    Throughout, $\catK(G) := \on{K}_b(\mathsf{perm}(kG)^\natural)$. We denote the tensor unit of $\catK(G)$ by $k := k[0]$ for shorthand (as opposed to $\mathbb{1}$) and the tensor product by $\otimes$ (as opposed to $\otimes_k$). Given a module $M$, $M[i] \in \catK(G)$ denotes the chain complex with $M$ in homological degree $i \in \Z$ and zero in all other degrees. We write $M := M[0]$ if the context is clear. The shift functor on $\catK(G)$ is also denoted $-[1]$. We denote by $\on{H}^\sbull(G) = \on{H}^\sbull(G;k)$ the $\N$-graded cohomology ring over $k$ in $\on{D}_b(kG) := \on{D}_b(\mathsf{mod}(kG))$ (including both odd and even-degree shifts for $p$ odd).

    \subsection*{Acknowledgments} The author is indebted to Paul Balmer and Martin Gallauer for numerous conversations and unwieldy email threads about their work, comments on earlier versions of this work, and introducing the author to the world of tensor-triangular geometry. He is especially thankful to Martin for poking holes in numerous previous iterations of the construction. He is grateful to Robert Boltje for his mentorship and suggestions during the course of the author's Ph.D., and the UC Santa Cruz Department of Mathematics for their generous graduate student travel support, which enabled numerous conversations leading to the ideas in this paper to happen. He thanks the referee for their exceptionally careful reading, numerous helpful suggestions, and for catching a few significant errors in prior versions. Finally, he thanks Cat \& Cloud Coffee of Santa Cruz, CA (on Swift St.) for their warm atmosphere, strong coffee, and impeccable breakfast burritos, countless of which were consumed during the lengthy inception of this project. At various stages of writing this paper, the author was partially funded by a UC Santa Cruz dissertation year fellowship, an AMS-Simons travel grant, and a US Junior Oberwolfach fellowship.

    During the research and writing stages of this paper, no generative AI tools were used. However, during a revision stage, the author used ChatGPT 6 Pro to double-check the correctness of arguments. The AI model discovered an critically false lemma (see \Cref{rmk:mistake}), which led to significant rewrites of latter sections by the author.

    \section{Preliminaries: Tensor-triangular geometry and permutation modules}

    For the first two sections, we review the essential prerequisites of the paper. Section 1 focuses more on tensor-triangular geometry and permutation modules, and Section 2 concerns endotrivial complexes and their topological connections.

    \subsection{Tensor-triangular geometry} We begin by reviewing the basic concepts of tensor-triangular geometry. For more detailed overviews of tensor-triangular geometry, we refer the reader to \cite{Bal05, Bal10b}.

    \begin{definition}
        Let $\catK$ be an essentially small tensor-triangulated category. That is, $\catK$ is an essentially small triangulated category with a symmetric monoidal structure such that for any $x \in \catK$, the functor $x \otimes -$ is an exact functor.
        \begin{enumerate}
            \item A triangulated subcategory $\catI \subseteq \catK$ is \emph{thick} if $\catI$ is closed under direct summands, and is a \emph{thick $\otimes$-ideal} if in addition, for all $x \in \catI$ and $y \in \catK$, $x \otimes y \in \catI$.
            \item A proper thick $\otimes$-ideal $\catP \subset \catK$ is \emph{prime} if for all $x,y \in \catK$, $x \otimes y \in \catP$ implies $x \in \catP$ or $y \in \catP$.
            \item Given a collection of objects $S \subseteq \catK$, we let $\langle S\rangle$ denote the thick $\otimes$-ideal generated by $S$. This is equivalently the intersection of all thick $\otimes$-ideals containing $S$.
        \end{enumerate}

        The following construction is due to Balmer \cite{Bal05}. Let $\Spc(\catK)$ denote the set of prime thick $\otimes$-ideals of $\catK$. For any $x \in \catK$, let $\supp(x) \subseteq \Spc(\catK)$ denote the subset \[\supp(x) := \{\catP \in \Spc(\catK) \mid x \not\in \catP\},\] the \textit{support} of $x \in \catK$. Given a collection of objects $S\subseteq \catK$, we set \[\supp(S) = \bigcup_{x\in S} \supp(x).\]
        The collection of all supports of objects induces a \emph{closed} basis on $\Spc(\catK)$, and a subset of $Z \subseteq \Spc(\catK)$ is of the form $Z = \supp(x)$ for some $x \in \catK$ if and only if $Z$ is closed with quasi-compact complement \cite[Proposition 2.14]{Bal05}. We call $\Spc(\catK)$ the \emph{Balmer spectrum} of $\catK$. It is well-known that $\Spc(-)$ is contravariantly functorial \cite[Proposition 3.6]{Bal05} and the \emph{universal support datum} for $\catK$ in a precise sense \cite[Theorem 3.2]{Bal05}. In particular, any \emph{tt-functor} (i.e., an exact $\otimes$-functor) induces a map on Balmer spectra.
    \end{definition}

    \begin{definition}
        Recall that given any $k$-linear category $\catC$, one can form its \emph{idempotent completion} $\catC^\natural$ (also called its \emph{Karoubi envelope}), which roughly speaking, closes the category $\catC$ under taking direct summands of objects. See for instance \cite{BS01} for triangulated categories. Given an open $U \subseteq \Spc(\catK)$ with closed complement $Z = \Spc(\catK) \setminus U$, we set \[\catK_Z := \{x \in \catK \mid \supp(x) \subseteq Z\},\] and define the localization \[\catK|_{U} :=\left(\catK/\catK_Z\right)^\natural.\] Then $\catK_Z$ is a thick $\otimes$-ideal of $\catK$, $\catK|_{U}$ is again a tensor-triangulated category, and we have a homeomorphism \[\Spc(\catK|_U) \cong U,\] see \cite[Proposition 1.11]{BF07}. This assignment endows $\Spc(\catK)$ with the structure of a locally ringed space \cite[Remark 6.4]{Bal05}.

        More generally, given any thick $\otimes$-ideal $\catI \subseteq \catK$ with canonical projection $\catK \mapsto \catK/\catI$, the induced map on Balmer spectra $\Spc(\catK/\catI) \to \Spc(\catK)$ induces a homeomorphism between $\Spc(\catK/\catI)$ and the subspace $\{\catP \in \Spc(\catK) \mid \catI \subseteq \catP\}$ of $\Spc(\catK)$ \cite[Proposition 3.11]{Bal05}.
    \end{definition}

    \begin{example}
        Let $G$ be a finite group (scheme), $k$ a field of characteristic $p$ prime, and $\on{D}_b(kG)$ denote the bounded derived category of $kG$-modules. We have a homeomorphism \[\Spc(\on{D}_b(kG)) \cong \Spec^h(\on{H}^\sbull(G;k)).\] This homeomorphism is realized by Balmer's \emph{comparison map} \cite{Ba10} \[\catP \mapsto \langle\zeta \in \on{H}^\sbull(G;k) \text{ homog.} \mid \cone(\zeta) \not\in \catP\rangle.\] We discuss an extended version of this comparison map in greater detail in the sequel. This theorem is due to Benson--Carlson--Rickard \cite{BCR97} for finite groups (in very different language than we present) and was extended to finite group schemes by Friedlander--Pevtsova \cite{FP07}. That the comparison map is a homeomorphism follows as a consequence of the classification of thick $\otimes$-ideals in both settings and Balmer's classification theorem \cite[Theorem 4.10]{Bal05}. Surjectivity follows from Noetherianity of the cohomology ring \cite{E61, FS97}, and a general tt-geometric fact \cite[Theorem 7.3]{Ba10}, but as far as we are aware, there is no known direct proof showing the comparison map is even injective.
    \end{example}

    \subsection{Permutation modules} Next, we review the basics of $p$-permutation modules. For a purely representation-theoretic overview of $p$-permutation modules and associated constructions, we refer the reader to \cite[Chapter 5]{L181} or \cite{CL23}.

    \begin{definition}
        A $kG$-module $M$ is:
        \begin{enumerate}
            \item a permutation module if $M$ admits a $G$-stable basis, or equivalently, $M \cong kX$ for some $G$-set $X$;
            \item a $p$-permutation module (where $p$ denotes the characteristic of $k$) if $M$ is a direct summand of a permutation module, or equivalently, $\Res^G_S M$ is a permutation module for some Sylow $p$-subgroup $S$ of $G$.

            In particular, if $G$ is a $p$-group, every $p$-permutation module is a permutation module. In this case, the indecomposable permutation modules are precisely those isomorphic to $kX$ for $X$ a transitive $G$-set. If $G$ is not a $p$-group, these modules may decompose further.
        \end{enumerate}
        We set $\catK(G) := \on{K}_b(\mathsf{perm}(kG)^\natural)$, the bounded homotopy category of $p$-permutation modules. The category $\catK(G)$ is a \emph{rigid} tensor-triangulated category, with duals given by $C^* := \Hom_k(C, k)$. Moreover, $\catK(G)$ is essentially small and is the compact part of a rigidly-compactly-generated (i.e., `big') tensor-triangulated category $\catD(G)$, the \emph{derived category of permutation $kG$-modules}. For this paper, $\catD(G)$ will not play much of a role.

        We note this category can be defined in the context of cohomological Mackey functors, Artin motives, or equivariant stable homotopy theory, and refer the reader to \cite{BG23b, F25} for details on these equivalences.
    \end{definition}

    There are two standard homomorphisms between permutation modules we will frequently consider.

    \begin{definition}
        Let $K \leq H \leq G$ be subgroups of $G$. The \textit{augmentation homomorphism} between the permutation modules $\on{aug}\colon  k[G/K] \to k[G/H]$ is the homomorphism induced by the $G$-set morphism given by $1K \mapsto 1H$. The standard example of the augmentation homomorphism is when $H = G$ and $K = 1$.

        There is an analogous \textit{coaugmentation homomorphism} $\on{coaug}\colon k[G/H] \to k[G/K]$ induced from the assignment $1H \mapsto \sum_{g \in [H/K]} gK$. This homomorphism of permutation modules is in no way induced from any $G$-set homomorphism.
    \end{definition}

    \begin{definition}
        By \cite{BG23}, we have a localization $\catK(G) \twoheadrightarrow \catK(G)/\catK_{ac}(G) \cong \on{D}_b(kG)$, where $\catK_{ac}(G)$ denotes the thick $\otimes$-ideal of $\catK(G)$ consisting of all acyclic complexes. We set $V_G := \Spc(\on{D}_b(kG))$, then \cite[Proposition 3.22]{BG25} asserts that this functor induces an open inclusion $V_G \hookrightarrow \Spc(\catK(G))$. The subset $V_G \subseteq \Spc(\catK(G))$ is called the \textit{cohomological open}. As noted in the introduction of \cite{BG25}, the `interesting' part of $\catK(G)$ is its closed complement, $\supp(\catK_{ac})$.
    \end{definition}

    \emph{Brauer quotients} are explicit constructions on $kG$-modules which mimic taking $H$-fixed points of $G$-sets when $H$ is a $p$-subgroup of $G$. Balmer--Gallauer adapt this construction to work on the level of big categories, dubbing them \textit{modular fixed points} functors.

    \begin{theorem}{\cite[Proposition 2.7]{BG25}}
        For every $p$-subgroup $H \leq G$ there exists a coproduct-preserving tt-functor on the big derived category of permutation modules \[\Psi^H = \Psi^H_G\colon  \catD(G) \to \catD(\Weyl{G}{H})\] such that $\Psi^H(kX) \cong k[X^H]$ for every $G$-set $X$. In particular, this functor preserves compacts and restricts to a tt-functor $\Psi^H\colon  \catK(G) \to \catK(\Weyl{G}{H})$.
    \end{theorem}

    The following conservativity theorem is an example of `local-to-global' phenomena in modular representation theory.

    \begin{theorem}{\cite[Theorem 2.11]{BG25}}\label{thm:conservative}
        The family of functors \[\{\catD(G) \xrightarrow{\Psi^H} \catD(\Weyl{G}{H}) \twoheadrightarrow \on{K}(\on{Inj}(k[\Weyl{G}{H}]))\}_{H \in \Sub_p(G)/G}\] indexed by conjugacy classes of $p$-subgroups $H \leq G$ is conservative, i.e., detects vanishing of objects. This restricts to a conservative family of functors  \[\{\catK(G) \xrightarrow{\Psi^H} \catK(\Weyl{G}{H}) \twoheadrightarrow \on{D}_b(k[\Weyl{G}{H}])\}_{H \in \Sub_p(G)/G}\] on the compact parts.
    \end{theorem}

    \begin{notation}
        Each tt-functor $\Psi^H$ induces a continuous map on spectra \[\Spc(\Psi^H) := \psi^H\colon  \Spc(\catK(\Weyl{G}{H})) \to \Spc(\catK(G)),\] and composing with the surjection $\check{\Psi}^H\colon  \catK(G) \twoheadrightarrow \on{D}_b(k[\Weyl{G}{H}])$ induces another continuous map on spectra \[\check{\psi}^H\colon  V_{\Weyl{G}{H}} \to \Spc(\catK(G)).\]
    \end{notation}

    These modular fixed point functors and the localization to the bounded derived category are the keys to deducing $\Spc(\catK(G))$ as a set. The Balmer spectrum, as a set, is built entirely from these $p$-local cohomological opens.

    \begin{theorem}{\cite[Theorem 2.10]{BG25}}
        Every point of $\Spc(\catK(G))$ is the image $\check\psi^H(\frakp)$ of a point $\frakp \in V_{\Weyl{G}{H}}$ for some $p$-subgroup $H \leq G$ unique up to $G$-conjugation, i.e., $\check\psi^H(\frakp) = \check\psi^{H'}(\frakp')$ if and only if there exists a $g \in G$ such that ${}^gH = H'$ and ${}^g\frakp = \frakp'$. In other words, \[\Spc(\catK(G)) = \bigsqcup_{H \in \Sub_p(G)/G} \Spc(\on{D}_b(k[\Weyl{G}{H}])) \cong \bigsqcup_{H \in \Sub_p(G)/G} \Spec^h(\on{H}^\sbull(\Weyl{G}{H};k)).\]
    \end{theorem}

    \begin{notation}\label{not:primes}
        We denote by $\catP_G(H, \frakp) \in \Spc(\catK(G))$ the image of the prime $\frakp \in V_{\Weyl{G}{H}}$ under $\check{\psi}^H: V_{\Weyl{G}{H}} \to \Spc(\catK(G))$. We omit the $G$ subscript if the ambient group is clear. By \cite[Theorem 2.10]{BG25}, every prime is, up to $G$-conjugation, uniquely expressible in this way.
    \end{notation}

    The subset $\{\catP(H, \frakp) \mid \frakp \in V_{\Weyl{G}{H}}\} \subseteq \Spc(\catK(G))$ is generally not open nor closed when $H \neq 1$, but contains the closed point $\catP(H,0) \in \Spc(\catK(G))$. Explicitly, one can have an inclusion of primes $\catP(H, \frakp) \subseteq  \catP(H', \frakq)$ where $H >_G H'$, see \cite[Corollary 7.13]{BG25}.

    \section{Preliminaries: Endotrivial complexes and Borel-Smith functions}

    We next review the classification of endotrivial complexes, following \cite{M24c}.

    \begin{definition}\label{rmk:hmarks}
        A chain complex $C \in \catK(G)$ is \textit{endotrivial} if it is an invertible object in $\catK(G)$, i.e., $C^* \otimes C \cong k$. The set $\Pic(\catK(G))$ of all endotrivials up to isomorphism forms a group under $\otimes$, and is called the \emph{Picard group} of $\catK(G)$.

        By the K\"unneth formula, for any $p$-subgroup $H\leq G$, $\Psi^H(C)$ has nonzero homology in exactly one homological degree, with that homology having $k$-dimension one. Let $h_C(H) \in \Z$ denote the unique integer for which  $\on{H}_{h_C(H)}(\Psi^H(C)) \neq 0$. The next proposition essentially reduces the study of endotrivials to considering superclass functions.
    \end{definition}

    \begin{prop}{\cite[Definition 3.5, Theorem 3.6]{M24a}}\label{prop:2.2}The following hold.
        \begin{enumerate}
            \item The function $h_C\colon  \Sub_p(G) \to \Z$ is a well-defined superclass function.
            \item The assignment $h\colon  \Pic(\catK(G)) \to \on{CF}(G,p)$ given by $C \mapsto h_C$ is a well-defined group homomorphism, with kernel isomorphic to $\Hom(G, k^\times)$.
        \end{enumerate}
    \end{prop}

    \begin{definition}
        Given an endotrivial $C$, we call $h_C$ the \textit{h-mark function} of $C$, and call the homomorphism $h$ the \textit{h-mark homomorphism}. Under the hood, $h$ is nothing more than a numerical avatar of the conservative family of functors $\{\check{\Psi}^H\}_{H \in \Sub_p(G)}$. Since the invertible objects of $\on{D}_b(kG)$ consists of shifts of $k$-dimension one $kG$-modules, $h_C$ simply tracks the image of $C$ in $\on{D}_b(k[\Weyl{G}{H}])$ where $H$ runs through all $p$-subgroups of $G$.
    \end{definition}

    By \Cref{prop:2.2}, to determine $\Pic(\catK(G))$, it effectively suffices to compute $\im(h)$.

    \begin{definition}
        A \textit{Borel-Smith function} is a superclass function $f \in \on{CF}(G,p)$ satisfying the following three conditions, which we call the Borel-Smith conditions.
        \begin{enumerate}
            \item If $p$ is odd, then for any subquotient $T/S$ of $G$ of order $p$, $f(T) \equiv f(S) \mod 2$.
            \item If $p = 2$, then for any sequence of subgroups $H \trianglelefteq K \trianglelefteq L \leq N_G(H)$, with $[K:H] = 2$, $f(K) \equiv f(H) \mod 2$ if $L/H$ is cyclic of order 4 and $f(K) \equiv f(H) \mod 4$ if $L/H$ is quaternion of order 8.
            \item For any elementary abelian subquotient $T/S$ of $G$ of rank 2, the equality \[f(S) - f(T) = \sum_{S < X < T} \big(f(X) - f(T)\big)\] holds.
        \end{enumerate}
        The collection of Borel-Smith functions forms an additive subgroup $\on{CF}_b(G,p)$ of $\on{CF}(G,p)$. In fact, if $G$ is a $p$-group, then under the identification $\on{CF}(G) \cong B^*(G) := \Hom(B(G), \Z)$, where $B(G)$ denotes the \emph{Burnside ring} of $G$, $\on{CF}_b(G)$ is a rational $p$-biset subfunctor of $\on{CF}(G)$. See \cite[Proposition 3.7]{BoYa07} for details.

    \end{definition}

    \begin{theorem} We have a complete characterization of endotrivials.
        \begin{enumerate}
            \item \cite[Theorem 4.6]{M24c} Let $G$ be a $p$-group. The h-mark homomorphism is an isomorphism onto the subgroup of Borel-Smith functions \[h\colon  \Pic(\catK(G)) \cong \on{CF}_b(G).\]
            \item \cite[Corollary 6.4]{M24c} Let $G$ be a finite group. The h-mark homomorphism has image $\on{CF}_b(G,p)$ and induces a split exact sequence \[0 \to \Hom(G,k^\times) \to \Pic(\catK(G)) \to \on{CF}_b(G,p) \to 0.\]
        \end{enumerate}
    \end{theorem}

    \begin{remark}
        The classification of endotrivials for a $p$-group can be deduced in multiple ways; one way is by identifying them with representation spheres, which we discuss in the following section. This method was considered by Bachmann in his dissertation \cite{Bac16}. Another way, the method we used in \cite{M24c}, is by recovering a short exact sequence constructed by Bouc--Yal\c{c}in \cite{BoYa07} relating Borel-Smith functions to the \emph{Dade group} of a finite group. This group parametrizes \emph{endopermutation} $kG$-modules, modules $M$ for which $M^* \otimes M$ is permutation. The classification of these modules for $p$-groups was completed by Bouc \cite{Bo06}.

        The classification of endotrivials for a non-$p$-group follows from some technical results concerning the homomorphism induced by restriction $\Res^G_H\colon \Pic(\catK(G)) \to \Pic(\catK(H))$ for subgroups $H$ containing a Sylow $p$-subgroup of $G$ deduced in \cite{M24b}.
    \end{remark}

    \begin{example}\label{ex:exampleendotriv}
        Let $G$ be a $p$-group. The first nontrivial example of an endotrivial, and those which make the magic happen in the elementary abelian case, are as follows (c.f. \cite[Definition 12.3]{BG25}). If $p = 2$, one has an endotrivial for $C_2$ given by $kC_2 \to k,$ with $k$ in homological degree 0 and the nonzero differential given by the augmentation homomorphism. If $p$ is odd, one has an endotrivial for $C_p$ by truncating a periodic resolution of the trivial $kC_p$-module $k$, $kC_p \to kC_p \to k$ with $k$ in homological degree 0. For any subgroup $N \triangleleft G$ of index $p$, inflation yields the following complex \[u_N:= k[G/N] \to k[G/N]\to k\] Here $p$ is assumed odd, if $p = 2$ then the degree 2 term is deleted. In the notation of Balmer--Gallauer, we set \[2' := \begin{cases} 1 & \text{if }p = 2; \\ 2 & \text{if } p > 2; \end{cases}\] hence $u_N$ has length $2'$. We keep the $u_N$ notation throughout the paper.

        A more secluded example for $p = 2$ and $G = D_{2^n}$ (here, we mean the dihedral group with $2^n$ elements, rather than the dihedral group acting on the regular polygon with $2^n$ edges) is as follows. Let $H_1, H_2$ be nonconjugate, non-central subgroups of order 2. This choice is unique up to conjugacy and reordering. Then the following complex is endotrivial: \[kD_{2^n} \to k[D_{2^n}/H_1] \oplus k[D_{2^n}/H_2] \to k,\] with all maps between indecomposable permutation modules being augmentation homomorphisms, $k$ in homological degree 0.
    \end{example}

    Finally, we note a structural fact about endotrivial complexes; the trivial $kG$-module $k$ only appears up to homotopy once.
    \begin{prop}\label{prop:onetrivial}
        Let $G$ be a $p$-group and $C \in \Pic(\catK(G))$ be an indecomposable endotrivial complex. Then $C_{h_C(G)}$ contains the trivial $kG$-module $k$ as a direct summand exactly once, and for all $j \neq h_C(G),$ $C_j$ does not contain $k$ as a direct summand.
    \end{prop}
    \begin{proof}
        Since $\Psi^G(C) \cong k[h_C(G)]$, $C_{h_C(G)}$ must contain a direct summand isomorphic to $k$, and \cite[Lemma 5.8]{M24a} with $\calY := \Sub_p(G)\setminus \Syl_p(G)$ proves the rest.

    \end{proof}

    \subsection{Borel-Smith functions and representation spheres}

    Borel-Smith functions initially arose not in the context of representation theory, but equivariant topology. To close out our preliminaries, we consider their origin: representation spheres, following \cite{tD87}.

    \begin{definition}\label{def:decreasing}
        We say that a superclass function $f$ is \textit{effective} if the following holds: if $K \leq H$ are subgroups of $G$, then $f(K) \geq f(H)$. That is, $f$ is monotonically decreasing with respect to the poset of subgroups of $G$. We say that an endotrivial $C$ is \textit{effective} if its corresponding h-mark function $h_C$ is effective.
    \end{definition}

    We first recall some important facts about Borel-Smith functions. Let $G$ be a group and $F$ a field. Given a $FG$-module $V$, the \textit{dimension function} associated to $V$ is the superclass function \[\dim\colon  H \mapsto \dim_F V^H.\] If $F$ has characteristic 0, this induces a group homomorphism $R_F(G) \to \on{CF}(G)$. When $F = \R$, it is common to write $\on{RO}(G)= R_\R(G)$, where $\on{RO}(G)$ denotes the Grothendieck ring of \emph{representation spheres}, i.e., 1-point compactifications of real representations (see \cite[Page 12]{tD87}). In this case, the dimension function is equivalently computing the dimension of the representation sphere $S^{V^H} = (S^V)^H$ of $V^H$.

    \begin{theorem}{\cite[Theorem 5.4 and Theorem 5.13, pages 211 and 216]{tD87}}\label{thm:imgofdimfunction}
        Let $G$ be a nilpotent group. The image of $\dim\colon  RO(G) \to \on{CF}(G)$ is the group of Borel-Smith functions $\on{CF}_b(G)$. Moreover, if $f$ is an effective Borel-Smith function returning non-negative values, there exists a real representation $V$ for which $\dim(V) = f$.
    \end{theorem}

    The following is immediate.

    \begin{corollary}
        Let $G$ be a nilpotent group and suppose $f$ is a Borel-Smith function. Then $f$ can be expressed as the difference of two effective Borel-Smith functions. In particular, every endotrivial is isomorphic in $\catK(G)$ to the product of an effective endotrivial and the dual of an effective endotrivial.
    \end{corollary}

    Therefore, it suffices to work with effective endotrivials, which we will see are particularly well-behaved. In fact, when $G$ is a nilpotent group, the dimension function gives a \textit{canonical} effective basis of $\on{CF}_b(G)$, which arises from the real irreducible representations of $G$. Therefore, if $G$ is a $p$-group, we have a corresponding canonical $\Z$-basis of $\Pic(\catK)$ consisting of effective endotrivials.

    \begin{theorem}
        Let $G$ be a $p$-group, and $V_1, \dots, V_n$ denote the irreducible real representations of $G$. Then the set of associated dimension functions $\{f_1,\dots, f_n\}$, after removing duplicates, forms a $\Z$-basis of $\on{CF}_b(G)$. In particular, there is an associated canonical $\Z$-basis of $\Pic(\catK(G))$.
    \end{theorem}
    \begin{proof}
        Since the image of the dimension homomorphism is precisely $\on{CF}_b(G)$, it suffices to show that the set $\calB := \{f_1, \dots, f_n\}$ is linearly independent after removing duplicates. \cite[Proposition 5.9, page 213]{tD87} asserts that $\ker(\dim)$ is generated by elements of the form $V - \psi^k(V)$, where $V$ is an irreducible real representation, $\psi^k$ is the $k$-th Adams operation, and $k$ is coprime to $|G|$. Therefore, the duplicates in $\calB$ arise from Adams-conjugates, hence any set of real irreducible representations with no Adams-conjugates will correspond to a linearly independent set of Borel-Smith functions.

    \end{proof}

    \begin{definition}\label{def:irreds}
        We call the Borel-Smith functions (resp. endotrivials) associated to the real irreducible representations of $G$ the \textit{irreducible} Borel-Smith functions (resp. endotrivials). These objects are necessarily effective.

        On the other hand, an \emph{indecomposable} endotrivial means what one would expect. It is an easy exercise to verify that in $\catK(G)$, for every element $C \in \Pic(\catK(G))$, there is a unique indecomposable representative of the isomorphism class of $C$.
    \end{definition}

    \begin{remark}
        Given a real representation $V$ of $G$ with corresponding character $\chi$, we have an equality \[\dim_\R V^H = \frac{1}{|H|}\sum_{h \in H} \chi(h)\] providing a computational method of computing the basis of $\on{CF}_b(G)$.

        Checking character tables of $G$ a $p$-group of normal $p$-rank one shows the basis of $\Pic(\catK(G))$ obtained in \cite[Section 6]{M24a} coincides with the canonical $\Z$-basis of $\on{CF}_b(G)$. We remark that this is rather surprising, as the computations of \cite{M24a} were performed entirely ad-hoc, as they occurred well before the classification of $\Pic(\catK(G))$.
    \end{remark}

    \begin{notation}
        Let $G$ be a $p$-group. Given a $\R G$-module $V$, let $f_V$ and $C_V$ denote the corresponding effective Borel-Smith function and effective indecomposable endotrivial respectively. That is, $f_V := \dim V$ and $C_V$ is the unique indecomposable endotrivial for which $h_{C_V} = \dim V$. \Cref{def:cwcplx} will explain how to obtain $C_V$ from $S^V$ directly.
    \end{notation}

    \section{Chain complexes over the orbit category}

    The connection between representation spheres, endotrivial complexes and Borel-Smith functions lets us discern that effective endotrivials have additional structure: they are induced from representation spheres via Bredon ($G$-equivariant) homology. $G$-CW-complexes produce chain complexes of free modules over the \textit{orbit category $\Gamma_G$}, which themselves induce particularly nice objects of $\catK(G)$, well-behaved with respect to local-global constructions. We first review the orbit category, following \cite{Y16}.

    \begin{definition}\label{def:orbitcat}
        Let $\Gamma_G$ denote the \textit{orbit category} of $G$. The objects of $\Gamma_G$ are transitive $G$-sets $G/H$ for subgroups $H\leq G$, and the morphisms from $G/H$ to $G/K$ are $G$-set homomorphisms $G/H \to G/K$. In particular, $\Hom_{\Gamma_G}(H,K)$ is empty unless a conjugate of $H$ is a subgroup of $K$. We may write a subgroup $H \leq G$ to denote the object $G/H \in \Gamma_G$ for shorthand.

        A (right) $k\Gamma_G$-module $M$ is a contravariant functor from the category $\Gamma_G$ to the category of $k$-modules. By identifying $\Aut_{\Gamma_G}(G/H)$ with $\Weyl{G}{H}$, $M(H)$ has $k[\Weyl{G}{H}]$-module structure. The category of finitely generated $k\Gamma_G$-modules, denoted $\catmod(k\Gamma_G)$ is abelian.

        Given a $G$-set $X$, we define the $k\Gamma_G$-module $k[X^?]$ as the module with value at $G/H$ given by the subspace $k[X^H]$, with the obvious induced maps. A module over $\Gamma_G$ is \textit{free} if it is isomorphic to a direct sum of modules of the form $k[(G/K)^?]$. Let $\catproj(k\Gamma_G)$ denote the full subcategory of $\catmod(k\Gamma_G)$ consisting of free $\Gamma_G$-modules (if $G$ is a $p$-group, then by the Yoneda lemma, every projective $\Gamma_G$-module is free, see \cite[Definition 3.1]{Y16}).
    \end{definition}

    The next propositions follow immediately from definition of $k\Gamma_G$-modules. We refer to these facts as the \textit{stabilizers grow} conditions. This property and name were first suggested to the author by Robert Boltje.

    \begin{prop}\label{prop:stabilizersgrow}
        Let $X$ and $Y$ be two $G$-sets. Any $k\Gamma_G$-homomorphism $f\colon k[X^?] \to k[Y^?]$ satisfies the following property: given any $x \in X$ and subgroup $H \leq G$, if $x \in X^H$, then $f(x) \in k[Y^H]$. In particular, if $X$ and $Y$ are transitive $G$-sets, then $\Hom_{\catmod(k\Gamma_G)}(k[X^?], k[Y^?])$ has $k$-basis induced from $\Hom_{G-\on{set}}(X,Y)$.
    \end{prop}

    \begin{prop}\label{cor:stabilizersgrow}
        Let $k[X^?]$ and $k[Y^?]$ be two free $k\Gamma_G$-modules and let $f\colon  k[X^?] \to k[Y^?]$ be a $k\Gamma_G$-module homomorphism. Let $K \leq H$ be two subgroups of $G$. For any $m \in k[X^H]$, if $f^H(m) = n \in k[Y^H]$, then regarding $m$ as an element of $k[X^K]\supseteq k[X^H]$ of $k[N_G(H)\cap N_G(K)]$-modules, $f^K(m) = n \in k[Y^K]$.
    \end{prop}

    \begin{remark}
        As a result, each free $k\Gamma_G$-module has a filtration by fixed points. Therefore, we may fix a `canonical permutation basis' in accordance with this stratification. We will make considerable use of this fact in the sequel.

        Given a $G$-set $X$ and two subgroups $K \leq H$ of $G$, there is a canonical inclusion homomorphism of $k[N_G(H) \cap N_G(K)]$-modules $i\colon k[X^H] \to k[X^K]$ associated to the $G$-set homomorphism $G/K \to G/H, 1K \mapsto 1H$. This comes associated with a (not necessarily unique) projection map $p\colon  k[X^K] \to k[X^H]$ satisfying $p\circ i = \id$ and $p(m) = 0$ if $m \in \Psi^H(k[X^K]) = 0$. After choosing a canonical permutation basis, we obtain an corresponding projection map associated to the basis.
    \end{remark}

    For any subgroup $H\leq G$, we have an obvious functor $-(H)\colon  \catproj(k\Gamma_G) \to \catperm(k[\Weyl{G}{H}])$ from the assignment $k[X^?] \to k[X^H]$. Modular fixed points behave as one would hope.

    \begin{prop}\label{prop:brauerquotientidentification}
        Let $P$ be a $p$-subgroup of $G$. Then the following diagram commutes up to natural isomorphism.
        \begin{figure}[H]
            \centering
            \begin{tikzcd}
                \catproj(k\Gamma_G) \ar[d, "-(1)"] \ar[dr, "-(P)"]\\
                \catperm(kG) \ar[r, "\Psi^P"] & \catperm(k[\Weyl{G}{P}])
            \end{tikzcd}
        \end{figure}
    \end{prop}
    \begin{proof}
        This follows from the natural isomorphism $k[X^P] \cong \Psi^P(kX)$.

    \end{proof}

    We now explain how to obtain an endotrivial complex from a representation sphere, via a construction topologists will be familiar with.

    \begin{definition}\label{def:cwcplx}
        Let $X$ be a $G$-CW-complex. The reduced cellular chain complex of $X$ over the orbit category is the functor $\tilde{C}(X^?) := \tilde{C}(X^?; k)$ from the orbit category $\Gamma_G$ to the category of chain complexes of $k$-modules. This gives rise to a chain complex of free $k\Gamma_G$-modules \[\tilde{C}(X^?) := \cdots \to C_i(X^?) \xrightarrow{d_i} C_{i-1}(X^?) \to\cdots \to C_0(X^?) \xrightarrow{\epsilon} \underline{k} \to 0,\] where $\underline{k}$ denotes the constant functor with values $k(H) = k $ for all subgroups $H\leq G$, and $\epsilon$ denotes the augmentation homomorphism. By convention $\underline{k}$ is in homological degree -1. This is essentially the construction of reduced Bredon homology with coefficients in $k$, see \cite{I73}. We denote chain complexes of free $k\Gamma_G$-modules by $C^?$, and their evaluations at the subgroup $H\leq G$ by $C^H$ to inspire an aura of fixed points.

        Given any chain complex of free $k\Gamma_G$-modules, evaluation at $1$ produces a chain complex of $p$-permutation $kG$-modules. This assignment induces a tt-functor \newline $(-)^1\colon\on{K}_b(\catproj(k\Gamma_G)) \to \catK(G)$, and the analogous commutative diagram to \Cref{prop:brauerquotientidentification} commutes.

        Let $G$ be a $p$-group. If $S^V$ is a representation sphere of $G$, it may be realized as a $G$-CW-complex by \cite{Ill83}, hence producing a chain complex of free $k\Gamma_G$-modules. Therefore for any subgroup $H \leq G$, $\tilde{C}(S^H)$ is the corresponding effective endotrivial complex of permutation $k[\Weyl{G}{H}]$-modules. It is clear there is an equality \[h_C(H) = \dim_\R V^H \] for all $H \leq G$.

        Conversely, by identifying an endotrivial $C$ with its h-mark homomorphism $h_C$, \Cref{thm:imgofdimfunction} asserts $C$ is up to shift homotopy equivalent to $\tilde{C}(S^1)$ for some representation sphere $S$ of $G$ with dimension homomorphism $\dim(S) = h_C$ by \Cref{prop:brauerquotientidentification}. Therefore, up to homotopy, every effective endotrivial arises from a chain complex of free $k\Gamma_G$-modules. We will show that one can remove ``up to homotopy'' in the sequel.
    \end{definition}

    \begin{construction}\label{cons:triangulation}
        Given a representation sphere $S^V$, we assume any $G$-CW complex structure arises from an equivariant triangulation of the sphere see e.g. \cite{Ill78}, for which stabilized simplices are pointwise stable. In particular, this uses that $S^V\cong_G S(V \oplus \R)$, where the latter term is the \emph{unit sphere} of $V \oplus \R$. So it suffices to give an equivariant triangulation of the unit sphere. For example, one may take the barycentric subdivision of the boundary of a $G$-invariant spanning polytope of the unit sphere.

        From this stipulation, the following is clear: if $V^K = V^H$ for two subgroups $H, K$ of $G$, then the induced triangulations, hence $G$-CW complex structures, on $(S^V)^H$ and $(S^V)^K$ are identical.
    \end{construction}

    \subsection{Removing contractible summands}\label{sec:removingsummands}

    We prove a couple technical results for $p$-groups that show we can remove contractible summands from chain complexes of free $k\Gamma_G$-modules if and only if we can for the corresponding chain complexes of permutation $kG$-modules (in particular, for effective endotrivials).

    \begin{remark}
        A chain complex of free $k\Gamma_G$ modules has in each homological degree a canonical permutation basis induced from the natural inclusions $X^H \hookrightarrow X^K$ for $K \leq H \leq G$, and each differential respects the stabilizers grow condition. Moreover, it is easy to see that given a chain complex $C$ of permutation $kG$-modules, if one can choose a canonical permutation basis of $G$-sets in each degree such that the stabilizers grow condition holds with respect to this basis, then $C$ may be realized as the image of a chain complex of free $k\Gamma_G$-modules evaluated at the subgroup $1 \leq G$. In this sense, the chain complexes of permutation $kG$-modules satisfying the stabilizers grow condition on differentials with respect to some chosen permutation basis are exactly those which arise from chain complexes of free $k\Gamma_G$-modules.

        If $C$ is a chain complex of permutation modules, we say $C$ \textit{arises from a chain complex of free $k\Gamma_G$-modules} if $C$ is in the essential image of the `evaluation at $G/1$' tt-functor $(-)^1\colon\on{K}_b(\catproj(k\Gamma_G)) \to \catK(G)$.

    \end{remark}

    \begin{lemma}\label{lem:nonkernelimpliessubiso}
        Let $G$ be a $p$-group and let $f\colon k[X^?] \to k[Y^?]$ be a free $k\Gamma_G$-homomorphism with $X$ transitive. Suppose the unique $k$-dimension 1 submodule $T \subseteq k[X^1]$ satisfies $T \not\subseteq \ker(f^1)$. Then there exists a direct summand $N$ of $k[Y^1]$ such that $p_N \circ f$ is an isomorphism.
    \end{lemma}
    \begin{proof}
        This follows from the standard fact that for indecomposable free $k\Gamma_G$-modules $k[(G/H)^?]$, $k[(G/K)^?]$, if there exists a (pointwise) injective morphism $f'\colon k[(G/H)^?]\to k[(G/K)^?]$, then $H,K$ are $G$-conjugates, and $f'$ is an isomorphism. Since $T \not\subseteq \ker(f^1)$, $f$ is (pointwise) injective, and projecting onto each direct summand of $k[Y^?]$, one such projection must be injective, hence an isomorphism.

    \end{proof}

    \begin{lemma}\label{prop:ctrctablesummandoforbitcc}
        Let $G$ be a $p$-group and suppose $C$ is a bounded chain complex of permutation $kG$-modules arising from a chain complex $C^?$ of free $k\Gamma_G$-modules. Let $i$ be an integer and suppose $T \subseteq C_i$ is a $k$-dimension one submodule. If $T \not\subseteq \ker(d_i)$, then there exists a nonzero contractible direct summand $K$ of $C$.

    \end{lemma}
    \begin{proof}
        Let $C_i = kX_1 \oplus \cdots \oplus kX_n$ denote the direct sum decompositions of $C_i$ given by the canonical permutation basis afforded by $C^?$. Then there exists an $j\in \{1,\dots, n\}$ for which, given any nonzero $t \in T$, we have $0 \neq (d_i\circ p_j)(t)$, where $p_j$ denotes projection onto $kX_j$. Then \Cref{lem:nonkernelimpliessubiso} implies that there exists an indecomposable direct summand $kY$ of $C_{i-1}$ for which $d_i$ restricted to $kX_j \to kY$ is an isomorphism. This forces the existence of a contractible direct summand of $C$ by a standard homological algebra argument, as desired.

    \end{proof}

    The next proposition demonstrates that if a chain complex $C^?$ of free $k\Gamma_G$-modules satisfies that $C^1$ has a contractible summand (as a chain complex of $kG$-modules), then $C^?$ itself has a contractible summand (as a chain complex of $k\Gamma_G$-modules). This fact will be rather useful for studying effective endotrivials. It is unclear if a non-technical proof of this fact exists, but we suspect that for chain complexes arising from $G$-CW-complexes, it may be so.

    \begin{theorem}\label{prop:removecontractiblesummands}
        Let $G$ be a $p$-group and suppose $C$ is a bounded chain complex of permutation $kG$-modules arising from a chain complex $C^?$ of free $k\Gamma_G$-modules. If $C$ contains a contractible direct summand, then there exists a direct sum decomposition $C = K \oplus D$ with $K$ contractible such that $D$ also arises from a chain complex $D^?$ of free $k\Gamma_G$-modules.

        In particular, $C$ has a contractible direct summand of permutation modules if and only if $C^?$ has a contractible direct summand.
    \end{theorem}
    \begin{proof}
        To prove this, it suffices to find a decomposition of $C$ into direct summands $K,D$ as stated and show that each homological degree of $D$ has a canonical permutation basis satisfying the stabilizers grow condition. If $C$ is indecomposable, there is nothing to show, so assume $C$ decomposes. First, since $C$ has a contractible direct summand, there exists an integer $i \in \Z$ and $k$-dimension one submodule $T \subseteq C_i$ such that $T \not\in \ker(d_i)$. Let $C_i = kX_1 \oplus \cdots \oplus kX_m$ and $C_{i-1} = kY_1\oplus \cdots \oplus kY_n$ be the direct sum decompositions into canonical bases, with each $X_j, Y_j$ a transitive permutation module. After projecting $T$ onto each $kX_j$, it follows that there exists an index $l$ for which $d_i|_{X_l}$ is injective. Choose an index $l$ such that the stabilizer of $X_l$ is minimal with respect to stabilizing group order.

        We modify the canonical bases of $C_i$ and $C_{i-1}$ as follows. First, \Cref{lem:nonkernelimpliessubiso} implies there exists a (non-unique) index $h \in \{1, \dots, n\}$ such that $p_{h}\circ d_i\circ i_l$ is an isomorphism, where $i_l$ denotes inclusion into $kX_l$ and $p_{h}$ denotes projection onto $kY_{h}$. We replace the basis elements of $Y_{h}$ with the basis elements $Y_{h}'=\{d_i(x) \mid x \in X_l\}$. Then $Y_h'$ is a $G$-set isomorphic to $X_l$, and since the projection of $kY_h'$ onto $kY_{h}$ is an isomorphism, this forms a new basis of $C_{i-1}$. Finally, $d_{i-1}(kY_h') = 0$ since $d_{i-1}\circ d_i = 0$, so the resulting basis still respects the stabilizers grow property. Set $D_{i-1} := kY_1 \oplus \cdots \oplus \widehat{kY}_h \oplus \cdots \oplus kY_n$ and $K_{i-1} := kY_h'$.

        Set $K_i := kX_l$, by construction $d_i|_{K_i}$ is an isomorphism and maps canonical permutation basis to canonical permutation basis. We next modify the canonical permutation basis of $C_n$ to obtain a direct summand $D_i$ for which $d_i(D_i) \subseteq D_{i-1}$. Set $C_{i}' := kX_1 \oplus \cdots \oplus \widehat{kX}_l \oplus \cdots \oplus kX_m$, so we have $C_i = C_i' \oplus K_i$. For each canonical permutation basis element $x$ of $C_i$ (i.e., some $x \in X_j$ for $j \neq l$) replace $x$ with $x' := x - ((d_i|_{kX_l})\inv\circ p_{h}'\circ d_i)(x)$, where $p_h'$ denotes projection onto $kY_h'$ and $(d_i|_{kX_l})\inv$ denotes the inverse of the isomorphic projection of $kX_l$ onto $kY_h'$. A straightforward computation shows $d_i(x') \in D_{i-1}$. By construction, $((d_i|_{kX_l})\inv\circ p_{h}'\circ d_i)(x)\in K_i$, and the collection $X'$ of the $x'$ forms a $G$-set, with $ X' \cong X_1 \sqcup \cdots \sqcup \widehat{X}_l \sqcup \cdots \sqcup X_m$. Set $D_i := kX'$, then we have a direct sum decomposition $C_i = D_i \oplus K_i$. Finally, for all $j \neq i, i-1$, set $D_j = C_j$. We have a direct sum decomposition $C = K \oplus D$. Finally, since all maps used to determine $x'$ are induced from $k\Gamma_G$-homomorphisms, the resulting differential $d_i$ satisfies the stabilizers grow property with respect to the new basis.


        It remains to show the stabilizers grow property holds for the differentials $d_{i-1}$ and $d_{i+1}$ with respect to the same canonical bases for $C_{i+1}$ and $C_{i-1}$. The condition holding for $d_{i-1}$ is straightforward, since the only modified basis elements in $C_{i-1}$ now belong to $\ker(d_{i-1})$. Similarly, $\im(d_{i+1}) \subseteq D_i$ since $K$ is contractible and lives in homological degrees $i$ and $i-1$. Under the decomposition $x' :=  x - ((d_i|_{kX_l})\inv\circ p_{h}'\circ d_i)(x)$ into $C_i'$ and $kY'_h$ respectively, it follows by construction that if a canonical permutation basis element $z \in C_{i+1}$ satisfies $d_{i+1}(z) = \sum a_j x_j + m$ with $m \in kX_l$ and each $x_j$ a canonical permutation basis element of $C_i'$, then we have $d_{i+1}(z) = \sum a_j x_j'$, where $x_j'$ is the refined basis element of $X'$ corresponding to $x_j$. Thus $d_{i+1}$ also satisfies the stabilizers grow condition, as desired.

    \end{proof}

    \begin{corollary}
        Let $G$ be a $p$-group. The tt-functor \[(-)^{1}: \on{K}_b(\catproj(k\Gamma_G)) \to \catK(G)\] induced by evaluation at the trivial subgroup $1$ of $G$ is conservative, i.e., detects vanishing of objects.
    \end{corollary}
    \begin{proof}
        \Cref{prop:removecontractiblesummands} implies the functor $(-)^{1}: \on{K}_b(\catproj(k\Gamma_G)) \to \on{K}_b(\catperm(kG))$ is conservative, and the canonical inclusion $\on{K}_b(\catperm(kG)) \to \on{K}_b(\catperm(kG))^\natural \cong \catK(G)$ is conservative as well.

    \end{proof}

    In particular, if $G$ is a $p$-group, then every indecomposable effective endotrivial arises from a free chain complex of $k\Gamma_G$-modules.

    \begin{theorem}\label{cor:indecendotrivsarefree}
        Let $G$ be a $p$-group and let $C$ be an indecomposable effective endotrivial of $kG$-modules. Then $C$ arises from a chain complex $C^?$ of free $k\Gamma_G$-modules. In particular, any effective endotrivial arises from a chain complex of free $k\Gamma_G$-modules.
    \end{theorem}
    \begin{proof}
        Recall every effective endotrivial (after a possible shift) corresponds to a real representation $V$ of $\R G$-modules. The representation sphere $S^V$ is a $G$-CW-complex, therefore produces a (reduced) chain complex $C^?$ for which $C^1$ is an endotrivial with the same h-marks as $C$. By \Cref{prop:removecontractiblesummands}, we can remove all contractible summands from $C^?$ until $C^1$ is indecomposable, and since two indecomposable endotrivials with the same h-marks are isomorphic, the result follows. The last statement follows easily by considering a direct sum decomposition of an effective endotrivial into its indecomposable endotrivial summand and contractible summand.

    \end{proof}

    We also prove a tightening of \Cref{prop:onetrivial} for effective endotrivials. The examples in \Cref{ex:exampleendotriv} demonstrate this, but it is not hard to find counterexamples of the proposition for $C$ non-effective.

    \begin{corollary}\label{prop:degzeroisk}
        Let $G$ be a $p$-group and $C \in \Pic(\catK(G))$ be an indecomposable effective endotrivial complex. Then $C_{h_C(G)} = k$.
    \end{corollary}
    \begin{proof}
        Without loss of generality we assume $h_C(G) = 0$. By \Cref{prop:onetrivial}, $C_0$ has a direct summand isomorphic to $k$ and for all $i \neq 0$, $k$ is not a direct summand of $C_i$. Furthermore, since $C$ is effective, by a standard $p$-permutation argument (see \cite[Proposition 5.8.11]{L181}), $C_i = 0$ for negative $i$. If $h_C \equiv 0$ then $C = k[0]$ so there is nothing to show. Otherwise, there exists a subgroup $H$ such that $h_C(H) > 0$, and in particular $h_C(1) > 0$.

        Now assume for contradiction that $C_0$ has another indecomposable summand $M \cong k[G/K]$ with $K < G$. Since $h_C(1) > 0$, $p_M\circ d_1$ is surjective, and because $C$ arises from a complex of free $k\Gamma_G$-modules, the fiber of $p_M \circ d_1$ in $C_1$ consists of indecomposable modules $k[G/K']$ with $K' \leq K$. Here $p_M$ denotes projection onto $M$.

        If $h_C(K) = 0$, then monotonicity implies $h_{\Psi^K(C)} \equiv 0$ so $\Psi^K(C) \cong k$. Therefore, there must be a summand $M'$ of $C_1$ also isomorphic to $k[G/K]$ in the fiber of $p_M \circ d_1$ to `cancel' $\Psi^K(M)$ via homotopy equivalence. But by \cite[Proposition 5.10.7]{L181} the composition $p_M \circ d_1 \circ i_{M'}$ is an isomorphism if and only if $\Psi^K(p_M \circ d_1 \circ i_{M'})$ is. Here $i_{M'}$ denotes the inclusion from $M'$ into $C_1$.  Therefore a contractible summand splits off (c.f. \cite[Lemma 9.2]{BM23}), contradicting indecomposability of $C$.

        Therefore, $h_C(K) > 0$ must hold. So $\Psi^K(d_1)$ is surjective, and there must also be a summand $M'$ of $C_1$ in the fiber of $p_M \circ d_1$ such that $\Psi^K(p_M\circ d_1\circ i_{M'})$ is surjective. Now a similar argument as before yields a contractible summand of $C$, a contradiction. We've exhausted all cases, thus $C_0 = k$.

    \end{proof}

    This final corollary will be critical in the sequel: it shows that nonzero homomorphisms $k \to C[s]$ for $C$ indecomposable and arising from a complex of free $k\Gamma_G$ modules, can be drawn without concern. If $C$ is endotrivial, \cite{Gal25} refers to such morphisms as \emph{global sections}.

    \begin{corollary}\label{cor:allthhekoms}
        Let $G$ be a $p$-group and let $C^?$ be a chain complex of free $k\Gamma_G$-modules with no contractible summands (e.g. an indecomposable effective endotrivial). For every integer $s \in \Z$ and submodule $T \subseteq C^1_s$ isomorphic to $k $, there exist a chain complex homomorphism $f\colon k  \to C[-s]$ with $\im(f_0) = T$. In particular, $\dim_k\Hom_{\catK(G)}(k, C^1[-s])$ is equal to the number of indecomposable direct summands of $C_s$.
    \end{corollary}
    \begin{proof}
        The existence of such a homomorphism $f$ is equivalent to the inclusion $T \subseteq \ker(d_s),$ and this inclusion holds since if not, \Cref{prop:ctrctablesummandoforbitcc} implies the existence of a contractible summand of $C^1$, hence a contractible summand of $C^?$, which cannot occur. The last statement is straightforward.

    \end{proof}

    \section{An open cover of the Balmer spectrum via endotrivials}

    \textbf{Notation:} \emph{For the rest of the paper,} we assume $G$ is a \textit{(finite) $p$-group}.

    \subsection{Constructing forerunners} We construct our `$p$-local quasi-isomorphisms,' the forerunners, for effective endotrivials.

    \begin{prop}\label{prop:rep_sphere_props_transferred}
        Let $G$ be a $p$-group, $V$ a real representation of $G$, and $C^? := \tilde{C}((S^V)^?)$ denote the associated cellular chain complex with differentials $d^?_i$. Set $H \leq G$ and $n := h_{C^1}(H)$. Then for all $i > n$, the permutation $kG$-module $C^1_i$ contains no indecomposable direct summands isomorphic to $k[G/K]$ with $K \geq H$. Consequently, the maximal nonzero homological degree of $\Psi^H(C^1) \cong C^H$ is $n$.
        Moreover, every indecomposable summand of $C^1_n$ isomorphic to $k[G/K]$ for any $K \geq H$ is isomorphic. We have $h_C(H) = h_C(K)$, and $K$ is maximal with respect to this property.
    \end{prop}
    \begin{proof}
        The first statement follows from the identification $S^{V^H} = (S^V)^H,$ and since $S^{V^H}$ is a $n$-sphere, it therefore contains no $i$-cells for $i > n$. Therefore $C^H_i = 0$, so the indecomposable direct summands of $C^1_i$ are killed by $\Psi^H$, and therefore are all transitive permutation modules isomorphic to $k[G/L]$ for some subgroup $L \not\geq_G H$ of $G$.

        To show that every indecomposable summand of $C^1_n$ with isotropy subgroup containing $H$ is in fact isomorphic to $k[G/K]$ for a fixed subgroup $K \geq H$ of $G$, suppose there is a direct summand of $C^1_n$ isomorphic to $k[G/K]$, with $K$ assumed to be maximal. By the previous statements, we must have that $h_{C^1}(H) = h_{C^1}(K)$. Indeed, $h_{C^1}$ is monotone, and $h_{C^1}(H) > h_{C^1}(K)$ would assert no summand of $C^1_n$ could be isomorphic to $k[G/K]$, a contradiction. Therefore, $V^K = V^H$ as vector spaces. Consequently, $(S^V)^H = (S^V)^K$ by \Cref{cons:triangulation}, so every $n$-cell that is $H$-fixed must also be $K$-fixed, so every indecomposable summand of $C^1_n$ not killed by $\Psi^H$ is also not killed by $\Psi^K$. By maximality of $K$, we are done.

    \end{proof}

    \begin{lemma}\label{prop:singletopmodule}
        Let $G$ be a $p$-group and let $C$ be an indecomposable effective endotrivial. Set $n := h_C(1)$. Then $n$ is the highest nonzero homological degree of $C$, $C_n$ is indecomposable, and if $V$ is a real representation of $G$ satisfying $\dim V \equiv h_C$ with kernel $N$, then $C_n \cong k[G/N]$.
    \end{lemma}
    \begin{proof}
        First, $V$ is equivalently a faithful $\R [G/N]$-module, so it suffices to assume $N = 1$ after replacing $G$ with $G/N$. From \Cref{prop:rep_sphere_props_transferred}, choosing $H = 1$ obtains that $C_n$ is the top homological degree of $C$. Moreover, since $V$ is faithful, it follows that there must be at least one $n$-cell that is acted upon faithfully by $G$, in which case \Cref{prop:rep_sphere_props_transferred} implies all such cells are acted upon faithfully. Consequently, $C_n$ is projective. It remains to show $C_n$ is indecomposable. Suppose for contradiction $C_n \cong P_1 \oplus P_2$ for projective $kG$-modules $P_1, P_2$. Since $\dim_k H_n(C) =1 $, either $d_n|_{P_1}$ or $d_n|_{P_2}$ is injective, hence split injective, so a contractible chain complex containing $P_1$ or $P_2$ splits off from $C$, contradicting indecomposability of $C$. Thus, $C_n$ is indecomposable projective, as desired.
    \end{proof}

    \begin{corollary}\label{cor:what_stabilizers_can_exist}
        Let $C$ be an indecomposable effective endotrivial complex, let $H\leq G$, and set $n := h_C(H)$. Then for all $i > n$, the permutation $kG$-module $C_i$ contains no indecomposable direct summands isomorphic to $k[G/K]$ with $K \geq_G H$. Therefore, the maximal homological degree of $\Psi^H(C)$ is $n$. Moreover, $\Psi^H(C)$, after removing contractible summands, contains a unique direct summand in degree $n$, and if $V$ is a real representation satisfying $\dim V \equiv h_C$, then this summand is a permutation module with stabilizer $\ker(V^H)$.
    \end{corollary}
    \begin{proof}
        This easily follows from \Cref{prop:singletopmodule} and \Cref{prop:rep_sphere_props_transferred}, using that every endotrivial contains a unique endotrivial direct summand, and all other summands are contractible.

    \end{proof}

    \begin{remark}\label{rmk:bgmorphisms}
        We recall the morphisms $a_N,b_N,c_N$ constructed by Balmer--Gallauer, \cite[Definition 12.3]{BG25}. Recall that given $N \triangleleft G$ a subgroup of $G$ of index $p$, one has the endotrivial of length $2'$ inflated from $C_p$ \[u_N = k[G/N] \to k[G/N] \to k,\] inflated from a truncated periodic resolution of the trivial module $k$.
        The above example is in the setting of $p > 2$. Then we have $2' + 1$ morphisms in $\catK(G)$, $a_N\colon k \to u_N$, $b_N\colon k \to u_N[-2']$, and if $p > 2$, $ c_N\colon k \to u_N[-1]$. The morphism $a_N$ in degree 0 is the identity map $k \xrightarrow{\id} k$. The morphisms $b_N$ and $c_N$ (if it exists) are induced by the coaugmentation homomorphism $k \to k[G/N]$.

        For all $H \leq G$, exactly one of $\Psi^H(a_N)$ or $\Psi^H(b_N)$ is a quasi-isomorphism, c.f. \cite[Proposition 12.9]{BG25}. On the other hand, $c_N$ is $\otimes$-nilpotent.
    \end{remark}

    We now generalize these maps for all endotrivials, as announced. This is an explicit construction which will be referenced later on.

    \begin{theorem}\label{thm:construction}
        Let $G$ be a $p$-group and let $C$ be an effective endotrivial. For every subgroup $H \leq G$, there exists a chain complex homomorphism \[\iota^H_C\colon k \to C[-h_C(H)]\] such that $\Psi^H(\iota^H_C)$ is a quasi-isomorphism.

        Moreover, there is an chain complex $C' \cong C$ arising from the cellularization of a representation sphere, and chain complex homomorphism $\iota^H_{C'}$ identifying with $\iota^H_C$ in $\catK(G)$, such that the image of $\iota^H_{C'}$ in $C'_{h_C(H)}$ is contained in a direct sum of indecomposable modules each isomorphic to $k[G/K]$, where $K \geq H$ is a subgroup of $G$ maximal with respect to the property that $h_C(H) = h_C(K)$.
    \end{theorem}

    \begin{proof}
        By \Cref{thm:imgofdimfunction}, it suffices to assume after degree shift that $C = C^1$, where $C^? = \tilde{C}((S^V)^?)$ for a representation sphere $S^V$ satisfying $\dim V \equiv h_C$. Denote the differentials of $C^?$ by $d^?$ and set $n := h_C(H)$. Denote the $G$-set of $n$-cells in (a triangulation of) $S^V$ by $X$, and let $X = X_1 \sqcup X_2$, where $X_1$ is all elements of $X$ with isotropy subgroup conjugate to a subgroup of $G$ containing $H$. Hence $(X_2)^H = \emptyset$ $(X_1)^H=X^H$, and with $C^1_n = kX = kX_1 \oplus kX_2$, we have $\Psi^H(kX_2) = 0$. By \Cref{prop:rep_sphere_props_transferred}, $X_1$ is moreover a sum of isomorphic transitive $G$-sets isomorphic to $k[G/K]$ for some subgroup $K \geq H$ of $G$, maximal with respect to the property that $h_C(H) = h_C(K)$.

        Now $C^H \cong \Psi^H(C)$ by \Cref{prop:brauerquotientidentification} and has maximal nonzero homological degree $n$, with $\ker(d^H_n) \cong k$. Moreover, $X^H$ consists of $n$-cells which have isotropy subgroup a conjugate to $K$ containing $H$.

        First, we consider the case when $K =H$. Then $\Weyl{G}{H}$ acts freely on $X^H$ and if $X_1$ is transitive, so is $X^H$ (see e.g. \cite[3.4]{Bro85}). We claim the following: \[\bigcup_{g\in [G/N_G(H)]} g\cdot X^H = X_1, \text{  and  } X^H \cap g\cdot X^H = \emptyset \text{   for all $g\in G\setminus N_G(H)$.}\] Since $X_1$ is a sum of isomorphic transitive $G$-sets, one reduces to the case of $X_1$ transitive, in which case $X_1 \cong G/H$ and $X^H \cong N_G(H)/H (= \Weyl{G}{H})$, and from here the claim follows.

        Define the $kN_G(H)$-module homomorphism $\tilde{\iota}^H_C$ by the following composition: \[k \cong \ker(d^H_n) \hookrightarrow C_n^H \hookrightarrow C_n^1.\] Note we also have an inclusion $ \ker(d^H_n)\subseteq \ker(d^1_n)$ from functoriality of $C^?\colon \Gamma_G\op \to \on{Ch}^b(\mathsf{mod}(k))$. We define $\iota^H_C\colon k \to C[-h_C(H)]$ in homological degree 0 as follows: \[\iota^H_C := \tr^G_{N_G(H)}(\tilde{\iota}^H_C) = \sum_{g\in [G/N_G(H)]} {}^g\tilde{\iota}^H_C\] and the zero map in all other degrees. We have from $G$-linearity that $\iota^H_C$ is a well-defined chain complex homomorphism, i.e., $\im((\iota^H_C)_0) \subseteq \ker(d_n^1)$. Moreover, that $\Psi^H(\iota^H_C)$ is a quasi-isomorphism follows from the fact that $X^H \cap g(X^H) = \emptyset$ for all $g \not\in N_G(H)$ and the identification $\Psi^H(C^1) \cong C^H$.

        Next, if the isotropy subgroup $K$ does not satisfy $K = H$, we have $K> H$. We define $\iota^H_C := \iota^K_C$, where $\iota^K_C$ is defined by the above (note the partition $X = X_1 \sqcup X_2$ with $K$ replacing $H$ does not change). We need to verify that $\Psi^H(\iota^K_C)$ is a quasi-isomorphism; however, this follows since $S^{V^K} = S^{V^H}$, therefore as complexes of $k$-modules, $\tilde{C}((S^V)^K) = \tilde{C}((S^V)^H)$. The final statement is clear from construction and the assumption that $C = C^1$.

    \end{proof}

    \begin{remark}\label{rmk:careiniota}
        We name the maps constructed in \Cref{thm:construction} \emph{forerunners}, and say $\iota^H_C$ is the \emph{$H$-forerunner of $C$}, as the maps herald the $p$-local isomorphisms $\check\Psi^H(C) \cong k[h_C(H)]$.

        \begin{enumerate}
            \item
            Whether a `canonical' construction of the forerunner homomorphism exists is a pertinent question. Such a construction may allow for better understanding of these morphisms.
            For instance, with our current construction it is unclear whether the following natural property holds: given two effective endotrivials $C_1, C_2$ and subgroup $H \leq G$, do we have an identification, up to scaling, between the images of $\iota^H_{C_1 \otimes C_2}$ and $\iota^H_{C_1}\iota^H_{C_2}$ in $\Hom_{\catK(G)}(k, C_1 \otimes C_2)?$

            By \Cref{prop:singletopmodule} and \Cref{prop:degzeroisk}, any indecomposable effective endotrivial $C$ satisfies that $C_{h_C(1)}$ is indecomposable and $C_{h_C(G)} = k$. It follows that $\iota^H_{C_1 \otimes C_2} = \iota^H_{C_1}\iota^H_{C_2}$ up to scaling for $H = 1$ and $H = G$, and consequently, for every effective endotrivial admitting at most 2 distinct h-marks.

            \item One has to take care in how the forerunners are constructed, as such maps that become quasi-isomorphisms locally are non-unique (not even up to choice of identification as in \cite[Remark 12.5]{BG25}). For instance, let $G = D_{2^n}$ for $n \geq 3$, and $C$ be the endotrivial of \Cref{ex:exampleendotriv}. \[C:= k D_{2^n} \xrightarrow{} k[ D_{2^n}/H_1] \oplus k[ D_{2^n}/H_2] \xrightarrow{} k ,\] where $H_1$ and $H_2$ are non-conjugate non-central subgroups of order 2 and the differentials are induced by augmentation homomorphisms. In this case, we could choose $\iota^{H_1}_C$ and $\iota^{H_2}_C$ to be the inclusion \[\iota\colon  k  \to k[D_{2^n}/H_1] \oplus k[D_{2^n}/H_2],\quad 1\mapsto  \left(\sum_{g\in [ D_{2^n}/H_1]} gH_1, \sum_{g\in [ D_{2^n}/H_2]} gH_2\right).\] In this case, both $\check{\Psi}^{H_1}(\iota)$ and $\check{\Psi}^{H_2}(\iota)$ are isomorphisms, a desired property. However, from the construction in \Cref{thm:construction}, one observes that $\iota^{H_1}_C$ and $\iota^{H_2}_C$ satisfy $\Psi^{H_1}(\iota^{H_2}_C) = 0$ and vice versa, since neither $H_1 \geq_G H_2$ or vice versa.

            More generally, one may observe from \Cref{thm:construction} that the forerunners associated to $C := C_{\R G}$ are all distinct. Indeed, one has $h_C(H) = [G:H]$, so in \Cref{thm:construction}, after choosing a suitable representative complex, every forerunners $\iota^H_C$ has image contained in a direct sum of $k[G/H]$'s. This observation plays a key role in constructing the open cover of $\Spc(\catK(G))$ stratifying the closed points.

            \item \label{rmk:conjugate_forerunners} By construction, given two subgroups $K,H$ of $G$, $\iota^H_C$ and $\iota^K_C$ may coincide.
            For instance, it is apparent from construction that $\iota^H_C = \iota^{{}^gH}_C$ for all $g \in G$., as long as the same associated cellularization is chosen. Alternatively, if $C = k$, then all forerunners coincide. We will assume for convenience that if $K,H$ are $G$-conjugate, then $\iota^H_C = \iota^K_C$ on the nose.

            \item Finally, a historical remark: in an earlier version of this paper, the forerunner construction was purely algebraic and only utilized the $k\Gamma_G$-module structure of an effective endotrivial, not the cellularization of a representation sphere. We posed the question of whether a topological construction, utilizing representation sphere structure, existed. The idea to consider a specific cellularization first appeared in \cite{GM26}, giving an answer to our question. However, an error was found in the original construction present in this paper, and we determined the easiest way to amend the issue was to use some of the ideas present in \cite{GM26}; this was done with permission of G\'omez. Balmer--Bouc--Dell'Ambrogio have developed another purely algebraic forerunner construction.
        \end{enumerate}

    \end{remark}

    The following lemmas give criteria for distinguishing forerunners.
    \begin{lemma}\label{prop:wheniotasareequal}
        Let $H, K$ be subgroups of $G$ and let $C$ be an effective endotrivial of $kG$-modules. If there exists a subgroup $B \leq G$ containing both $H,K$ up to conjugacy for which $ h_C(B) = h_C(H) = h_C(K)$, then $\iota^H_C = \iota^K_C$ (up to a possible scalar).
    \end{lemma}
    \begin{proof}
        This follows from \Cref{thm:construction} after choosing the associated representation sphere $S^V$ associated to $C$ and assuming $C = C_V$. If such a $B$ exists, then equality of h-marks gives after conjugating $H,K$ into $B$ that $V^H = V^B = V^K$ so the maximal isotropy subgroups of \Cref{thm:construction} agree and the forerunners coincide up to a scalar.

    \end{proof}

    \begin{lemma}\label{lem:samesizesubs}
            Let $H, K$ be non-conjugate subgroups of $G$ with same order and let $C$ be an effective endotrivial. Suppose $h_C(H) = h_C(K)$, and for every $g \in G$, we have $h_C(\langle H,{}^gK\rangle) < h_C(H)$. Then $\Psi^K(\iota^H_C) = 0$ and $\Psi^H(\iota^K_C) = 0$.
        \end{lemma}
    \begin{proof}
            This can be observed from \Cref{thm:construction} and \Cref{prop:rep_sphere_props_transferred}.
            Choose a representative $C$ arising from the cellularization of a representation sphere. Then the images of $\iota^H_C$ and $\iota^K_C$ are contained in a sum of transitive permutation modules isomorphic to $k[G/H']$ and $k[G/K']$ respectively, where $H' \geq H$ and $K'\geq K$ are determined up to conjugacy. Moreover, we have $C_{h_C(H)} = kX_1 \oplus kX_2 = kY_1 \oplus kY_2$, where $X_1$ is a sum of $k[G/H']$'s, $Y_1$, is a sum of $k[G/K']$'s, $X_2^H = 0$, and $Y_2^K = 0$.

            It follows that $\Psi^K(\iota^H_C) \neq 0$ if and only if one can choose $H' = K'$, and similarly for $\Psi^H(\iota^K_C)$. In particular, if this holds, then $H'$ contains both $H$ and $K$ up to conjugacy. However, for any $g \in G$, we have that $\langle H,{}^gK\rangle$ satisfies $h_C(\langle H,{}^gK\rangle) < h_C(H) = h_C(K)$, therefore $C_{h_C(H)}$ contains no transitive permutation modules with isotropy subgroups containing $\langle H,{}^gK\rangle$ by \Cref{prop:rep_sphere_props_transferred}. It follows that $H', K'$ cannot contain $\langle H,{}^gK\rangle$ for any $g$, hence $H'$ does not contain $K$ up to conjugacy, and similarly $K'$ does not contain $H$ up to conjugacy. The result follows.
    \end{proof}
%

    \subsection{An open cover of the spectrum} Using the forerunners, we construct the promised open cover of $\Spc(\catK(G))$. Recall that every prime of $\Spc(\catK(G))$ is of the form $\catP(H, \frakp) := \check\psi^H(\frakp)$ for some subgroup $H \leq G$ and $\frakp \in V_{\Weyl{G}{H}} := \Spc(\on{D}_b(k[\Weyl{G}{H}]))$, and $\catP(H, \frakp) = \catP(H', \frakp')$ if and only if $H' = {}^gH$ and $\frakp' = {}^g\frakp$ for some $g \in G$.

    \begin{construction}\label{con:open}
        With the notation of \cite[Definition 12.3]{BG25}, given any normal subgroup $N \leq G$ of index $p$, the maps $a_N$ and $b_N$ are examples of such morphisms (however, $c_N$ is not!). The endotrivial \[u_N = k[G/N] \to k[G/N] \to k \] (assuming $p$ odd) satisfies $a_N = \iota_{u_N}^H$ for any $H \not\leq N$ and $b_N = \iota_{u_N}^H$ for any $H \leq N$.

        Given this observation, we generalize the open cover of $\Spc(\catK(G))$ presented in \cite[Proposition 13.11]{BG25} as follows, with one key amendment. Let $\calB(G)$ denote the subset of the canonical $\Z$-basis of $\on{CF}_b(G)$ not induced from the trivial $\R G$-module $\R$, i.e., the effective endotrivials arising from irreducible real representation spheres, excluding $k[1]$. For any subgroup $H \leq G$, we define an open subset of $\Spc(\catK(G))$ via \[U(H) := \bigcap_{C \in  \calB(G)} \open(\iota_C^H),\] and \[U(H)' := U(H) \cap \open(\iota^H_{C_{\R G}}).\] Here, \[\open(f) := \open(\cone(f)) = \{\catP \mid f \text{ is invertible in } \catK(G)/\catP\}.\]

        Since $\iota^H_C = \iota^{{}^gH}_C$ for any $g \in G$, $U(H) = U({}^gH)$. Moreover, $\open(\iota_{k[1]}^H) = \Spc(\catK(G))$ for any subgroup $H \leq G$. Therefore, there is no harm in excluding $k[1] = C_{\R}$ from $\calB(G)$; this exclusion is mostly for notation and indexing purposes.
    \end{construction}

    Let us recall the closed points of $\Spc(\catK(G))$, i.e., the minimal primes with respect to inclusion. \cite[Corollary 7.31]{BG25} asserts that the minimal primes of $\Spc(\catK(G))$ are those $\catM(H) := \catP(H,0)$ for some subgroup $H$, where $0$ denotes the unique closed point of $V_{\Weyl{G}{H}}$. The closed point $\catM(H)$ is precisely the kernel of the residue tt-functor $\F^H := \Res_1^{\Weyl{G}{H}}\circ\Psi^H = \Psi^H \circ \Res^G_H$, see \cite[Definition 7.26]{BG25}, and $\catM(H)$ is also the kernel of $\check\Psi^H$.

    \begin{prop}\label{prop:closedptsinopen}
        Let $H$ be a $p$-subgroup of $G$ and $C$ be an effective endotrivial. Recall the residue tt-functor \[\mathbb{F}^H\colon=  \Res^{\Weyl{G}{H}}_1 \circ \Psi^H\colon  \catK(G) \to \on{D}_b(k)\] of the closed point $\catM(H)$. We have that $\mathbb{F}^H(\iota^H_C)$ is an isomorphism. Moreover, if $K$ is a subgroup of $G$ satisfying $h_C(K)\neq h_C(H)$, then $\mathbb{F}^K(\iota^H_C)$ is not an isomorphism.

        In particular, the opens $U(H)$ and $U(H)'$ contain $\catM(H)$, therefore the sets of opens $\{U(H)\}_{H \in \Sub_p(G)/G}$ and $\{U(H)'\}_{H \in \Sub_p(G)/G}$ are open covers of $\Spc(\catK(G))$.
    \end{prop}
    \begin{proof}
        The fact that $\F^H(\iota^H_C)$ is an isomorphism follows immediately since $\Psi^H(\iota^H_C)$ is a quasi-isomorphism by construction. Moreover, if $h_C(K) \neq h_C(H)$, $\mathbb{F}^K(\iota^H_C)$ cannot be an isomorphism since $k[h_C(H)]$ is not isomorphic to $\Psi^K(C)$ in $\on{D}_b(k[\Weyl{G}{K}])$.

        By conservativity, for any effective endotrivial $C$, $\iota^H_C$ is an isomorphism in $\catK/\catM(H)$, as $\F^H$ is the residue functor of $\catM(H)$, so $U(H)$ contains $\catM(H)$. Therefore by general tt-geometry,  $\{U(H)\}_{H \leq G}$ is an open cover of $\Spc(\catK(G))$, as every prime specializes to some $\catM(H)$, which are precisely the closed points of $\Spc(\catK(G))$.

    \end{proof}

    The open covers witness the fact that every endotrivial is a tt-line bundle.
    \begin{theorem}\label{cor:linebundle}
        Every element $C \in \Pic(\catK(G))$ is a tt-line bundle under the open cover $\{U(H)\}_{H \in \Sub_p(G)/G}$. In particular, for every subgroup $H \leq G$, we have an isomorphism $k\cong C[-h_C(H)]$ in the localization $\catK(G)|_{U(H)}$. Consequently, the same holds for $\{U(H)'\}_{\Sub_p(G)/G}$
    \end{theorem}
    \begin{proof}
        By construction of $U(H)$ and the forerunners, every irreducible endotrivial corresponding to an element of the canonical $\Z$-basis of $\on{CF}_b(G)$ is isomorphic to $k[h_C(H)]$ in $\catK(G)|_{U(H)}$ (after shifting the forerunner $h_C(H)$ homological degrees). As the set of irreducible endotrivials is a $\Z$-basis of $\Pic(\catK(G))$, every endotrivial $C \in \Pic(\catK(G))$ satisfies $C[-h_C(H)] \cong k$ in $\catK(G)|_{U(H)}$. The same holds for each $U(H)'$ since $U(H)' \subseteq U(H)$, so we have localizations $\catK(G) \twoheadrightarrow \catK(G)|_{U(H)} \twoheadrightarrow \catK(G)|_{U(H)'}$.

    \end{proof}

    \begin{remark}
        Given an arbitrary endotrivial $C$, \Cref{cor:linebundle} does not imply the existence of an $H$-forerunner, i.e., a homomorphism $f: k \to C[-h_C(H)]$ which descends to an isomorphism in $\catK(G)|_{U(H)}$. Such homomorphisms may not exist in general for non-effective endotrivials. For example, let $p = 2 $ and $G = C_2$, then the two-term endotrivial $C := k \xrightarrow{\on{coaug}} kC_2$ with $k$ in homological degree 0 satisfies $\Hom_{\catK(G)}(k, C) = 0$, but we have a local isomorphism $k[0] \cong C$ in the localization $\catK(G)|_{U(C_2)}$.
    \end{remark}

    In fact, the opens $\{U(H)'\}$ stratify the closed points.
    \begin{theorem}\label{cor:onlyclosedpoint}
        The prime $\catM(H)$ is the only closed point of $\Spc(\catK(G))$ contained in the open $U(H)'$.
    \end{theorem}
    \begin{proof}
        We previously showed in \Cref{prop:closedptsinopen} that $\catM(H) \in U(H)'$. Since $U(H)'$ is defined by iterating over the irreducible endotrivials and $C_{\R G}$, it suffices to show $\F^K(\iota^H_{C_{\R G}}) = 0$ for all nonconjugate pairs of subgroups $H,K$. If $h_C(H) \neq h_C(K)$ then this is shown in \Cref{prop:closedptsinopen}, so assume $h_C(H) = h_C(K)$.
        Here, we have $h_{C_{\R G}}(H) = [G:H]$, therefore from \Cref{thm:construction}, every subgroup $H$ is locally maximal with respect to all subgroups $K \geq H$ satisfying $h_C(H) = h_C(K)$. Therefore, every forerunner $\iota^H_{C_{\R G}}$ has image contained in a direct sum of permutation modules isomorphic to $k[G/H]$. From this, it is clear that $\Psi^H(\iota^K_{C_{\R G}}) = 0$, hence $\F^H(\iota^K_{C_{\R G}}) = 0$, and the result follows.

    \end{proof}

    \begin{remark}
        We additionally warn that \Cref{cor:onlyclosedpoint} does not imply that the open $U(H)'$ has a unique closed point, as a subspace of $\Spc(\catK(G))$. The open $U(H)'$ has a minimal closed point if and only if $U(H)'$ consists of all points specializing to $\catM(H)$. In general this may not occur; this can already be observed for elementary abelian $p$-groups.
    \end{remark}

    \begin{remark}\label{rem:nonconjex}
        We define our opens $U(H)'$ by intersecting with the forerunner associated to $C_{\R G}$ entirely to prove \Cref{cor:onlyclosedpoint}. At present, it is unclear to us whether the theorem holds without this additional intersection, i.e., for the opens $U(H)$. If one was able to show the property $\iota_{C_1}^H \otimes \iota_{C_2}^H = \iota_{C_1 \otimes C_2}^H$ held for every subgroup $H$ and effective endotrivials $C_1, C_2$, then the intersection would indeed be redundant, but we are not sure if this holds.

        We briefly remark on what can possibly go wrong; it turns out that in most cases, \Cref{cor:onlyclosedpoint} holds for $\{U(H)\}$ as well, but there is a class of possible exceptions. Let $K, H$ be a pair of non-conjugate subgroups of $G$. Say $K$ and $H$ are \textit{indistinguishable} if for all Borel-Smith functions $f$, $f(K) = f(H)$. If no such indistinguishable pairs exist in $G$, say \textit{non-conjugacy is detected in $G$.} If $K$ and $H$ are indistinguishable, we have $|K| = |H|$, since $\dim(\R G)^H = [G:H]$, the index of $H$ in $G$.

        Many $p$-groups have non-conjugacy detected, such as abelian groups trivially, all $p$-groups of normal $p$-rank one (this is computed indirectly in \cite[Section 6]{M24a}), and all groups of order at most $p^3$. However, indistinguishable pairs of subgroups exist. For instance, set $G := C_8 \rtimes (C_2 \times C_2)$, with generators $a,b,c$ satisfying \[a^8 = b^2 = c^2 = 1, \quad {}^ba = a\inv, \quad {}^ca = a^3, \quad bc=cb.\]
        This group has a GAP implementation of \texttt{SmallGroup(32,43)}; it is the \emph{holomorph} of $C_8$, i.e., a group of the form $H \rtimes \Aut(H)$. In particular, $C_2 \times C_2\cong \Aut(C_8)$ acts on $C_8$ faithfully. The group $G$ has two nonconjugate subgroups isomorphic to $V_4 = C_2 \times C_2$, \[H := \{1, b, c, bc\}, \quad K := \{1, a^2b, a^2c, bc\},\] whose individual elements are all conjugate: \[b\sim_G a^2b, \quad c\sim_G a^2c.\] Therefore, given any real representation $V$ of $G$ with character $\catP_V$, we have $\dim V^H = \dim V^K$, since \[\dim V^H = \frac{1}{|H|}\sum_{h \in H}\chi_V(h).\] Since $\im(\dim) = \on{CF}_b(G)$, every Borel-Smith function $f$ satisfies $f(H) = f(K)$. We thank Math.SE user \texttt{testaccount} for finding and sharing this example \cite{SEt}.

        For these groups, it is not necessarily clear that the targets of forerunners of irreducible endotrivials are enough to stratify the closed points of $\Spc(\catK(G))$. One needs the existence of an effective endotrivial $C$ such that $\F^H(\iota_C^K) = 0$ for indistinguishable subgroups $H,K$ of $G$, and irreducible endotrivials need not accomplish this feat.
    \end{remark}

    \subsection{The open $U(1)$} We turn to the open $U(1)$, which as it turns out, is precisely the cohomological open. Note that since $\iota_{C_1}^1 \otimes \iota_{C_2}^1 = \iota_{C_1 \otimes C_2}^1$ by \Cref{prop:singletopmodule} (see \Cref{rmk:careiniota}), $U(1) = U(1)'$.

    \begin{definition}
        We recall Balmer--Gallauer's \emph{Koszul complexes}: given a subgroup $H \leq G$, set \[\on{kos}_G(H) := \on{Ten}^G_H (0 \to k \xrightarrow{\id} k \to 0),\] where $\on{Ten}^G_H$ denotes \emph{tensor induction} of chain complexes, see \cite[Section 4.1]{Be982}. Note that while tensor induction is monoidal, it does not preserve contractibility, and hence is not functorial for $\catK(G)$. One has that the tt-ideal $\ker(\Res^G_H)$ is generated by $\on{kos}_G(H)$ by \cite[Proposition 3.21]{BG25}. In particular, $\supp(\on{kos}_G(H)) = \supp(\ker(\Res^G_H))$.
    \end{definition}

    \begin{remark}\label{rmk:needallinvs}
        \Cref{cor:onlyclosedpoint} does not hold for the open cover of \cite[Section 13]{BG25} for non-elementary abelian $p$-groups, demonstrating that to stratify $\Spc(\catK(G))$, one needs to consider more endotrivials than just the collection of $u_N$'s. For instance, \cite[Proposition 13.14]{BG25} states the closed complement of the open $U'(1)$ (writing $U'(1)$ to denote the open $U_G(1)$ of \cite[Proposition 13.14]{BG25}, in order to distinguish the open covers apart) is the support of $\on{kos}_G(F)$, where $F$ denotes the Frattini subgroup of $G$ (i.e., the intersection of all maximal subgroups of $G$). By \cite[Corollary 7.17]{BG25}, we have $\catM(H) \in \supp(\on{kos}_G(F))$ if and only if $H \not\leq_G F$, hence $\catM(H) \in U'(1)$ if and only if $H \leq_G F$. Therefore, $U'(1)$ contains a unique closed point of $\Spc(\catK(G))$ if and only if $F = 1$ if and only if $G$ is an elementary abelian $p$-group.
    \end{remark}

    Note that for any effective endotrivial $C$, $\iota^1_C$ is unique up to scalar since $C_{h_C(1)}$ is indecomposable by \Cref{prop:singletopmodule}.
    First, we need a lemma about \textit{faithful} endotrivials, i.e., endotrivials that arise from a faithful real representation $V$. In this case, we recover a generalization of \cite[Proposition 13.14]{BG25}.

    \begin{prop}\label{lem:faithfulendotriv}
        Let $C$ be an effective endotrivial arising from a faithful real representation $V$. Then $\cone(\iota^1_C)$ generates $\catK_{ac}(G)$ as a tt-ideal, and we have an equality \[\supp(\cone(\iota^1_C)) = \supp(\ker(\Res^G_1)) = \supp(\catK_{ac}(G)) = \supp(\on{kos}_G(1)).\]
    \end{prop}
    \begin{proof}
        We have that $\ker(V) = 1$, so by \Cref{prop:singletopmodule} if $n$ is the highest nonzero homological degree of $C$, $C_n \cong kG$. Now consider the dual complex $\cone(\iota^1_C)^*$. After shifting, this complex satisfies $C_1 \cong kG$ and $C_0 \cong k $. Therefore, we are in the situation of \cite[Corollary 3.20]{BG25}, so $\cone(\iota^1_C)^*$, hence $\cone(\iota^1_C)$, generates $\catK_{ac}(G)$ as a tt-ideal, and $\supp(\catK_{ac}(G)) = \supp(\cone(\iota^1_C))$, as desired.

    \end{proof}

    \begin{remark}\label{rem:kernels}
        Analogously, if $V$ is a real representation with kernel $N \trianglelefteq G$, then it is a faithful real representation of $G/N$, and the associated endotrivial $C_V$ is a faithful endotrivial complex of $k[G/N]$-modules.
        In this case, $\cone(\iota^1_{C_V}) = \cone(\iota^N_{C_V})$ (as they are the same map) and \[\supp(\ker(\Res^G_N)) =  \supp(\on{kos}_G(N)) = \supp(\Inf^G_{G/N}( \on{kos}_{G/N}(N/N))),\] with the rightmost equality implied by \cite[Lemma 3.17]{BG25}. Therefore we have \[\supp(\cone(\iota_{C_V}^N)) = \supp(\ker(\Res^G_N)),\] generalizing \cite[Lemma 13.2(b)]{BG25} to arbitrary endotrivials.
    \end{remark}

    \begin{theorem}\label{prop:openisac}
        Let $G$ be a $p$-group. The closed complement of the open $U(1)$ is the support of $\on{kos}_G(1)$, i.e., the closed support of the tt-ideal $\catK_{ac}(G) = \ker(\Res^G_1)$. In particular, $U(1)$ is equal to the cohomological open $V_G = \Spc(\on{D}_b(kG)).$
    \end{theorem}
    \begin{proof}
        Since $\calB(G) \cup \{k[1]\}$ is a basis of $\Pic(\catK(G))$, for any $\catP \in U(1)$, all endotrivials are isomorphic to a shift of the tensor unit in $\catK(G)/\catP$. In particular, there exists a faithful real representation $V$ (for instance the regular representation $\R G$) and a corresponding endotrivial $C_V$, for which $\iota_{C_V}^1\colon  k [h_{C_V}(1)] \to C_V$ is an isomorphism in $\catK(G)/\catP$ (in fact, an isomorphism in $\catK(G)|_{U(1)}$ by \Cref{cor:linebundle}). Therefore, $\open(\cone(\iota^1_{C_V})) \supseteq U(1)$. By \Cref{lem:faithfulendotriv}, we have $\supp(\on{kos}_G(1)) = \supp(\cone(\iota_{C_V}^1)) \subseteq U(1)^c$.

        Conversely, by definition we have \[U(1) = \bigcap_{C \in \calB(G)}\open(\iota^1_C).\]
        Therefore by \Cref{rem:kernels},
        \begin{align*}
            U(1)^c &= \bigcup_{C \in \calB(G)}\supp(\cone(\iota^1_C))\\
             &=\bigcup_{C \in \calB(G)} \supp(\on{kos}_G(\ker(C))),
        \end{align*}
        where $\ker(C)$ denotes the kernel of the real representation corresponding to $C$. For any subgroup $H \leq G$, one has by \cite[Corollary 7.17]{BG25} $\supp(\on{kos}_G(H)) \subseteq \supp(\on{kos}_G(1))$, thus $U(1)^c \subseteq  \supp(\on{kos}_G(1))$, thus equality holds.

    \end{proof}

    The following corollary is immediate by general tt-geometry, since $\catK(G)$ is rigid.

    \begin{corollary}
        We have a tt-equivalence $\catK(G)|_{U(1)} \cong \on{D}_b(kG)$.
    \end{corollary}

    In particular, the open $U(1)$ has a unique closed point as a subspace of $\Spc(\catK(G))$. However, we stress that the other opens $U(H)$ need not have unique closed points. In particular, the spaces $\Spc(\catK(G)/\catM(H)) \subseteq \Spc(\catK(G)|_{U(H)}) = U(H)$ need not be the same, as the former has a unique closed point, but the latter may not.

    \section{The remixed twisted cohomology ring}

    We continue to assume $G$ is a $p$-group. In this section, we introduce the twisted cohomology ring $\tH(G)$ for $G$.

    \begin{definition}\label{def:twistedcohomology}
        Let $ \N^{\calB(G)} = \{q\colon  \calB(G) \to \N\}$ be the \textit{monoid of twists}, i.e., tuples of non-negative integers indexed by elements of $\calB(G)$. Equivalently, the monoid of twists identifies with the submonoid of $\Pic(\catK(G))$ generated by $\calB(G)$. Consider the $(\Z\times \N^{\calB(G)})$-graded ring \[\tH(G) = \tH(G;k) := \bigoplus_{s\in \Z}\bigoplus_{q \in \N^{\calB(G)}} \Hom_{\catK(G)} \left(k , \bigotimes_{C \in \calB(G)} C^{\otimes q(C)}[s]\right).\] Multiplication is induced by the tensor product in $\catK(G)$. Note that only \textit{non-positive} shifts $s \leq 0$ produce non-zero homomorphisms. We call $\tH(G)$ the \textit{(permutation) twisted cohomology ring of $G$}, as in \cite{BG25}. It is convenient to write \[k(q) := \bigotimes_{C \in \calB(G)} C^{q(C)}\] for every twist $q \in \N^{\calB(G)}$ and thus abbreviate $\on{H}^{s,q}(G) = \Hom_{\catK(G)}(k , k(q)[s])$.

        This ring is commutative for $p = 2$ and graded-commutative for $p$ odd (see \cite[Remark 12.8]{BG25}) - note that only the shift $s$ plays a role in the graded commutativity, and not the twist $k(q)$.
    \end{definition}

    \begin{remark}\hfill
        \begin{enumerate}
            \item When $G$ is a finite elementary abelian $p$-group, we recover the twisted cohomology ring constructed in \cite{BG25}, since in this case the set of endotrivials $\{u_N\}_{N \in \calN}$ as defined in \cite[Example 12.1]{BG25} and the shift of the tensor unit $k[1]$ form precisely the canonical $\Z$-basis of $\Pic(\catK(G))$, as deduced in \cite[Section 6]{M24a}.
            \item This construction may not immediately adapt for non-$p$-groups, since in this setting it is not known if one has a canonical $\Z$-basis for $\Pic(\catK(G))$. A possible replacement we propose could be the subgroup of $\catK(G)$ of endotrivials arising from representation spheres for $G$. It is unknown if the subgroup of $\Pic(\catK(G))$ spanned by such endotrivials even has maximal rank; we expect these questions to be nontrivial.
        \end{enumerate}

    \end{remark}

    To examine the $k$-vector space $\on{H}^{s,q}(G)$, it suffices to choose an indecomposable representative of $k(q)$. This proposition also classifies all the \emph{global sections} of an effective endotrivial, in the sense of \cite{Gal25}.

    \begin{prop}\label{prop:identifyingkhoms}
        Let $G$ be a $p$-group, $C$ an effective endotrivial, and $s$ an integer. Let $f\colon  k  \to C[-s]$ be a nonzero chain complex homomorphism. Then $[f] \in \Hom_{\catK(G)}(k , C[-s])$ is null homotopic if and only if $\im(f_0) \subset C_s$ is contained in a contractible direct summand of $C$. In particular, if $C$ is an indecomposable effective endotrivial, then \[\on{H}^{-s,q}(G) := \Hom_{\catK(G)}(k , C[-s]) = \Hom_{\on{Ch}(\mathsf{perm}(kG)^\natural)}(k , C[-s]) = \Hom_{\mathsf{mod}(kG)}(k, C_{s}).\]
    \end{prop}
    \begin{proof}
        The converse implication is straightforward. Suppose $f \in \Hom_{\catK(G)}(k , C[-s]) = 0$, then there exists a homotopy $h\colon  k  \to C_{s+1}$ such that $f = d_{s+1} \circ h$. Choose a canonical basis of $C$, and consider the projection of the image of $f$ onto each indecomposable summand. Necessarily on each summand $k[G/H]$ whose projection is nonzero, the image must be the unique minimal nonzero submodule $k  \subseteq k[G/H]$. Since $f = d_{s+1} \circ h$, there exists an indecomposable direct summand $M$ of $C_{s+1}$ such that $p_{k[G/H]}\circ f = d_{s+1}\circ p_{M} \circ h$. Since the stabilizers grow condition holds for $C$, the only possible homomorphisms $M \to k[G/H]$ with $k  \subset M$ not in the kernel are isomorphisms, hence $M \cong k[G/H]$ and we may split off the contractible summand. An induction argument demonstrates $\im(f)$ is contained in a contractible direct summand of $C$ as desired. The rightmost equality follows from \Cref{cor:allthhekoms}.

    \end{proof}

    As in \cite{Ba10, BG25}, we have a canonical comparison map, the heart of the matter. Its existence follows identically to \cite[Proposition 13.4]{BG25}. See also \cite[Appendix A]{Gal25} for a much more general treatment of twisted cohomology rings.

    \begin{prop}\label{prop:comparisonmap}
        There is a continuous \textit{comparison map}
        \[\on{comp}_G\colon  \Spc(\catK(G)) \to \Spec^h(\tH(G))\] mapping a tt-prime $\catP$ to the ideal generated by the homogeneous $f \in \tH(G)$ whose cone does not belong to $\catP$. It is characterized by the fact that for all $f$, \[\on{comp}_G\inv(Z(f)) = \on{supp}(\cone(f)) = \{\catP \mid f \text{ is not invertible in } \catK(G)/\catP\}\] where $Z(f) = \{\mathfrak{p}\mid f \in \mathfrak{p}\}$ is the \emph{closed} subset of $\Spec^h(\tH(G))$ defined by $f$.
    \end{prop}

    \begin{remark}
        Recall for a morphism $f$ in $\catK(G)$ we write \[\open(f) := \open(\on{cone}(f)) = \{\catP \mid f \text{ is invertible in } \catK(G)/\catP\}.\] The open $\open(f)$ is the preimage of $\on{comp}_G$ of the principal quasi-compact open $Z(f)^c = \{\mathfrak{p} \mid f \not\in \mathfrak{p}\}$, and is equivalently the open locus of $\Spc(\catK(G))$ where $f$ is invertible.
    \end{remark}

    \subsection{Twisted cohomology under localization and tt-functors}

    We extend the results of \cite[Section 14]{BG25} regarding localization and functoriality. The following constructions rely on the forerunners constructed in \Cref{thm:construction}.

    \begin{definition}\label{def:cohlocalization}
        Let $H$ be a subgroup of $G$. Let $S_H \subset \tH(G)$ be the multiplicative subset generated by all forerunners $\iota^H_C$, where $C$ runs over the \emph{irreducible} endotrivials $C \in \calB(G)$. We define a $\Z$-graded ring \[\calO^\sbull_G(H) := \big(\tH(G)[S_H\inv]\big)_{\text{0-twist}}\] as the twist-zero part of the localization of $\tH(G)$ with respect to $S_H$. Explicitly, the homogeneous elements of $\calO^\sbull_G(H)$ consist of fractions $\frac{f}{g}$ where $f, g \in \tH(G)$ with the same $\calB(G)$-twist $q$, and $g$ is a product of `$H$-local isomorphisms' $g'\colon  k  \to C[s]$ with $C \in \calB(G)$. Thus, $\calO^\sbull_G(H)$ is $\Z$-graded by the shift only. The homological degree of $\frac{f}{g}$ is the difference $s-t$ between the shifts of $f$ and $g$.
    \end{definition}

    \begin{construction}\hfill\label{cons:localizations}
        \begin{enumerate}
            \item We can perform the central localization (c.f. \cite{Ba10}) of the whole category $\catK(G)$ \[ \calL_G(H) := \catK(G)[S_H\inv].\] In fact, this localization has idempotent completion $\catK(G)|_{U(H)}$ by the same argument as in \cite[Construction 14.12]{BG25} or \cite[Theorem 3.6]{Ba10}.

            Explicitly, the category $\calL_G(H)$ is the Verdier quotient of $\catK(G)$ by the tt-ideal $\langle \{\cone(g) \mid g \in S_H\}\rangle$. It has the same objects as $\catK(G)$ and morphisms $x \to y$ of the form $\frac{f}{g}$ where $g\colon  k  \to C[s]$ belongs to $S_H$ for $C$ an effective endotrivial and $f\colon  x \to C[s] \otimes y$ is any morphism in $\catK(G)$ with the same twist $C$ as the denominator $g$. Moreover, the $\Z$-graded endomorphism ring $\End^\sbull_{\calL_G(H)}(k )$ of the unit in $\calL_G(H)$ is the $\Z$-graded ring $\calO^\sbull_G(H)$.

            In particular, $\calO^\sbull_G(1) \cong \on{H}^\sbull(G)$ by \Cref{prop:openisac}. We refer the reader to \cite[Section 3]{Ba10} for this construction in abstract.

            \item Twisted cohomology $\tH(G)$ is graded over the monoid $\Z \times \N^{\calB(G)}$. Ring homomorphisms induced by tt-functors or localization will be homogeneous with respect to a certain homomorphism $\gamma$ on the corresponding grading monoids.

            Let $H \leq G$ be a subgroup and consider the central localization \[(-)_{U(H)}\colon  \catK(G) \twoheadrightarrow \calL_G(H).\] Here, the $H$-forerunners all become isomorphisms, yielding a homomorphism on the grading \[\gamma_{U(H)}\colon  \Z \times \N^{\calB(G)} \to \Z\] defined by $\gamma_{U(H)}(s,q) = s + h_{k(q)}(H)$ and one obtains a ring homomorphism \[(-)_{U(H)}\colon  \tH(G) \to \End^{\sbull}_{\calL_G(H)}(k) = \calO^\sbull_G(H),\] homogeneous with respect to the monoid homomorphism $\gamma_{U(H)}$.

            \item \label{cons:modularfixedptsontc}Let $H \trianglelefteq G$ be a normal subgroup. The tt-functor $\Psi^H\colon  \catK(G) \to \catK(G/H)$ maps effective endotrivials with h-mark at $G$ equal to 0 to effective endotrivials at $G/H$ equal to 0. In fact, since deflation preserves irreducibility of a representation, $\Psi^H$ induces a surjective map $\calB(G) \to \calB(G/H) \sqcup \{k\}$.

            This defines a homomorphism of graded monoids \[\gamma = \gamma_{\Psi^H}\colon  \Z \times \N^{\calB(G)} \to \Z \times \N^{\calB(G/H)}\] given by $\gamma(s,q) = (s, \overline{q})$ where $\overline{q}$ is given by the surjection $\N^{\calB(G)} \twoheadrightarrow \N^{\calB(G/H)}$ along the inclusion $\calB(G/H) \hookrightarrow \calB(G)$ induced from inflation. Therefore, modular fixed points defines a ring homomorphism $\Psi^H\colon \tH(G) \to \tH(G/H)$ homogeneous with respect to $\gamma_{\Psi^H}$. Because $\id_{\catK(G/H)} = \Psi^H \circ \Inf^G_{G/H}$ general algebraic geometry implies $\Spec^h(\tH(G/H))$ is homeomorphically embedded into $\Spec^h(\tH(G))$ via the map induced by $\Psi^H$ on spectra.

            \item Given a group homomorphism $\alpha\colon G' \to G$, restriction along $\alpha$ also defines a tt-functor $\alpha^*: \catK(G) \to \catK(G')$. This again defines a corresponding ring homomorphism $\alpha^*\colon \tH(G) \to \tH(G')$ homogeneous with respect to $\gamma_{\alpha^*}$.
            If $\alpha$ is surjective, then $\alpha^*$ is simply inflation $\Inf^G_{G/N}$ for $N = \ker(\alpha)$. If $\alpha$ is injective, then $\alpha^*$ is restriction. We write \[\rho^G_H := \left(\Res^G_H\right)^*\colon \Spec^h(\tH(H)) \to \Spec^h(\tH(G)).\] These are essentially the only two cases one needs to consider, as every group homomorphism factors as a surjection composed with an injection.
        \end{enumerate}
    \end{construction}

    \begin{remark}\label{ex:resnotclosed}
        For any subgroup $H \leq G$, the tt-functor $\Res^G_H\colon \catK(G) \to \catK(H)$ induces a closed map on Balmer spectra $\Spc(\catK(H)) \to \Spc(\catK(G))$ by \cite[Proposition 4.7]{BG25}; this follows by monadicity of the restriction/induction adjunction, see \cite{Ba16}. The same does not hold for the induced map on twisted cohomology $\rho^G_H\colon\Spec^h(\tH(H)) \to \Spec^h(\tH(G))$. Indeed, let $G = C_2$, then $\tH(C_2) \cong k[a,b]$ with $a,b$ the usual $a_1, b_1$ morphisms in $\catK(C_2)$:
        \begin{figure}[H]
            \centering
            \begin{tikzcd}
                k \ar[d, "b"] & & & k\ar[d, "a"]\\
                kC_2 \ar[r] & k & kC_2\ar[r] & k
            \end{tikzcd}
        \end{figure}
        Now consider the restriction $\Res^{C_2}_1\colon \tH(C_2) \to \tH(1)$. The image of the induced map on spectra $\rho^{C_2}_1\colon \Spec^h(\tH(1)) \to \Spec^h(\tH(C_2))$ is the point associated to the prime $\langle a\rangle \subset \tH(C_2)$, which is neither an open nor closed point of $\Spec^h(\tH(C_2))$.
    \end{remark}

    \begin{remark}\label{rmk:commsquares}
        Suppose $F\colon  \catK(G) \to \catK(G')$ is a tt-functor and the induced homomorphism $F\colon  \tH(G) \to \tH(G')$ is homogeneous with respect to $\gamma = \gamma_F\colon  \Z \times \N^{\calB(G)} \to \Z \times \N^{\calB(G')}$ (for instance, modular fixed points or restriction). Then the following square commutes, since $F(\cone(f)) = \cone(F(f))$.

        \begin{figure}[H]
            \centering
            \begin{tikzcd}
                \Spc(\catK(G')) \ar[r, "F^*"] \ar[d, "\on{comp}_{G'}"'] & \Spc(\catK(G)) \ar[d, "\on{comp}_{G}"]\\
                \Spec^h(\tH(G')) \ar[r, "F^*"] & \Spec^h(\tH(G))
            \end{tikzcd}
        \end{figure}

        Similarly, the localization functor $(-)|_{U(H)}$ obtains the following commutative square.

        \begin{figure}[H]
            \centering
            \begin{tikzcd}
                U(H) = \Spc(\calL_G(H)) \ar[r, hookrightarrow] \ar[d, "\on{comp}_{\calL_G(H)}"'] & \Spc(\catK(G)) \ar[d, "\on{comp}_{G}"]\\
                \Spec^h(\calO^\sbull_G(H)) \ar[r, hookrightarrow] & \Spec^h(\tH(G))
            \end{tikzcd}
        \end{figure}
        Here, the left hand vertical map is the usual graded comparison map from \cite{Ba10} for the tt-category $\calL_G(H)$.
    \end{remark}

    We have an analogous version of \cite[Proposition 14.21]{BG25}.

    \begin{prop}\label{prop:diagonalinjective}
        Let $H \leq G$ be a subgroup. Then we have a commutative square
        \begin{figure}[H]
            \centering
            \begin{tikzcd}
                V_{\Weyl{G}{H}} = \Spc(\on{D}_b(k[\Weyl{G}{H}])) \ar[r, "\check{\psi^H}"] \ar[d, "\on{comp}_{\on{D}_b(k[\Weyl{G}{H}])}"'] & \Spc(\catK(G)) \ar[d, "\on{comp}_G"]\\
                \Spec^h(\on{H}^\sbull(\Weyl{G}{H})) \ar[r] & \Spec^h(\tH(G))
            \end{tikzcd}
        \end{figure}
        If $H \trianglelefteq G$ is normal, then the diagonal is injective.
    \end{prop}
    \begin{proof}
        The above diagram factors as follows:
        \begin{figure}[H]
            \centering
            \begin{tikzcd}
                \Spc(\on{D}_b(k[\Weyl{G}{H}])) \ar[r, hookrightarrow] \ar[d, "\comp_{\on{D}_b(k[\Weyl{G}{H}])}"'] & \ar[r, "\psi^H"]\Spc(\catK(\Weyl{G}{H})) \ar[d, "\comp_{\Weyl{G}{H}}"] & \Spc(\catK(G))  \ar[d, "\comp_G"]\\
                \Spec^h(\on{H}^\sbull(\Weyl{G}{H})) \ar[r, hookrightarrow] & \Spec^h(\tH(\Weyl{G}{H})) \ar[r, "\psi^H"] & \Spec^h(\tH(G))
            \end{tikzcd}
        \end{figure}

        Here, the left-hand square arises from the central localization $(-)|_{U(H/H)}$ which by \Cref{prop:openisac} is identified with the localization $\catK(\Weyl{G}{H}) \to \on{D}_b(k[\Weyl{G}{H}])$, and the right-hand square is induced from the tt-functor $\Psi^H\colon \catK(G) \to \catK(\Weyl{G}{H})$. \Cref{rmk:commsquares} asserts both diagrams commute, thus the square commutes as well.

        If $H \trianglelefteq G$, then $\Psi^H\colon \tH(G) \to \tH(G/H)$ is surjective, as noted in \Cref{cons:localizations}(c), hence $\psi^H: \Spec^h(\tH(G/H)) \hookrightarrow \Spec^h(\tH(G))$ is injective. Additionally, the comparison map $\comp_{\on{D}_b(k[\Weyl{G}{H}])}$ is always a homeomorphism, thus the bottom left composition is injective, and therefore so is the diagram.

    \end{proof}

    \subsection{An identification for the image of $\comp_G$}

    We may consider the commutative diagram of \Cref{prop:diagonalinjective} for all conjugacy classes of subgroups of $G$ at once. Recall the conservative functor of \Cref{thm:conservative}, the collection of modular fixed points functors \[\{\check{\Psi}^H\}_{H \in \Sub_p(G)/G}\colon \catK(G) \to \prod_{H \in \Sub_p(G)/G} \on{D}_b(k[\Weyl{G}{H}]).\] For shorthand, we denote this functor by $\check\Psi$. Denote the cohomology ring of $\prod_{H \in \Sub_p(G)/G} \on{D}_b(k[\Weyl{G}{H}])$ by \[\on{dH}^{\sbull}(G) = \on{dH}^{\sbull}(G,k) = \prod_{H \in \Sub_p(G)/G} \on{H}^\sbull(\Weyl{G}{H},k).\] We have an induced map on cohomology from \Cref{cons:localizations}, also denoted by $\check\Psi$. \[\check\Psi\colon \tH(G) \xrightarrow{\prod \Psi^H} \prod_{H \in \Sub_p(G)/G} \tH(\Weyl{G}{H}) \to \prod_{H \in \Sub_p(G)/G} \calO_{\Weyl{G}{H}}^\sbull(H/H) = \on{dH}^\sbull(G).\] We denote the associated map on homogeneous ring spectra by $\check\psi^h\colon \Spec^h(\on{dH}^\sbull(G)) \to \Spec^h(\tH(G))$.
    The following observation is easy but crucial, as it allows us to understand the comparison map $\comp_G$ purely in terms of induced maps on (Zariski) spectra. Injectivity of the diagonal is the crux of \Cref{sec:inj}.

    \begin{prop}\label{prop:dbcommdiagram}
        The following diagram in $\on{Top}$ commutes.
        \begin{figure}[H]
            \centering
            \begin{tikzcd}
                \bigsqcup_{H \in \Sub_p(G)/G} V_{\Weyl{G}{H}}\ar[d, "\cong"] \ar[r, hookrightarrow, "\check{\psi}"] & \Spc(\catK(G)) \ar[d, "\on{comp}_G"]\\
                \Spec^h(\on{dH}^{\sbull}(G)) \ar[r, "\check\psi^h "] & \Spec^h(\tH(G))
            \end{tikzcd}
        \end{figure}
        In particular, the image of the comparison map $\on{comp}_G$ is as a set precisely the image of $\check\psi^h $.
    \end{prop}
    \begin{proof}
        Commutativity follows directly from commutativity of the diagrams of \Cref{prop:diagonalinjective} for all subgroups $H \leq G$. The latter statement follows since the top and left arrows are bijections of sets.

    \end{proof}

    \begin{remark}\label{rem:equivimage}
        Unwinding the commutative diagram in \Cref{prop:dbcommdiagram}, the image of the prime $\catP(H, \frakp) \in \Spc(\catK(G))$, with $\frakp$ regarded as an element of\\ $\Spec^h(\on{H}^\sbull(\Weyl{G}{H}))\cong V_{\Weyl{G}{H}}$, under the comparison map $\comp_G$ is \[\langle f \in \tH(G) \text{ homog.} \mid \check\Psi^H(f) \in \frakp\rangle.\]
    \end{remark}

    While the comparison map is in general not surjective, it has dense image.

    \begin{corollary}\label{cor:denseimage}
        The image of $\comp_G$ is dense in $\Spec^h(\tH(G))$.
    \end{corollary}
    \begin{proof}
        Let $f \in \on{H}^{s,q}(G)$ be a homogeneous element of $\tH(G)$. By \cite[Theorems 7.1, 7.2]{BG25}, if $\check\Psi(f) = 0$, then $f$ is $\otimes$-nilpotent. We claim this fact is sufficient to prove that $\im(\comp_G) = \im(\check\psi^h)$ is dense in $\Spec^h(\tH(G))$. Indeed, this follows by a standard argument. Recall that $\Spec^h(\tH(G))$ has an open basis given by supports $\supp(f)$ of \emph{homogeneous} $f\in \tH(G)$. Therefore, for all homogeneous $f \in \tH(G)$, if $\supp(f) \neq \emptyset$, then $f \not\in\frakp$ for some homogeneous $\frakp \in \Spec^h(\tH(G))$. Hence $f$ is not $\otimes$-nilpotent, and it follows that neither is $\check\Psi(f)$, so there exists a prime $\frakq \in \Spec^h(\on{dH}^\sbull(G))$ with $\check\Psi(f) \not\in\frakq$. Therefore $f \not\in \check\psi^h(\frakq)$, so $\check\psi^h(\frakq) \in \supp(f)$. Thus $\im(\check\psi^h)$ intersects every nonempty basic open set $\supp(f)$, as desired.

    \end{proof}

    \begin{remark}
        Note that the previous proof does not imply that every element of $\ker(\check\Psi)$ is $\otimes$-nilpotent. Indeed, as an example, let $G = C_2$, then $\tH(C_2)\cong k[a,b]$ as in \Cref{ex:resnotclosed}. We have $\Psi^1(b) = \id$ and $\Psi^{C_2}(b) = 0$, therefore the non-homogeneous element $b(b - \id) \in \tH(C_2)$ satisfies $b(b-\id) \in \ker(\check\Psi)$. However $b(b-\id)$ cannot be $\otimes$-nilpotent, as $\tH(C_2)$ contains no nilpotent elements. In particular, $\ker(\check\Psi)$ may be a non-homogeneous ideal of $\tH(G)$.
    \end{remark}

    \section{Injectivity of the comparison map}\label{sec:inj}

    We continue to assume $G$ is a $p$-group.
    The crux of injectivity of the comparison map is showing the diagonal of \Cref{prop:diagonalinjective} is injective when a subgroup $H \leq G$ isn't normal. For this next proposition, we introduce subscripts for modular fixed points to avoid confusion.

    \begin{prop}\label{prop:newhardlemma}
        Let $G$ be a $p$-group and $H \leq G$. For every homogeneous $f \in \tH(N_G(H))$, there exists a homogeneous $\hat{f} \in \tH(G)$ such that \[\check\Psi_{N_G(H)}^H(f) = \check\Psi_G^H(\hat{f})\in \on{H}^\sbull(N_G(H)/H).\]
    \end{prop}
    \begin{proof}
        Let $f \colon k \to D[s]$ be the homogeneous element, with $D$ an effective endotrivial. First, we replace $f$ with a homogeneous $f' \colon k \to D'[s']$ such that $D'$ is the restriction of an effective endotrivial of $kG$-modules and that $\check\Psi_{N_G(H)}^H(f') = \check\Psi_{N_G(H)}^H(f)$. Since $D$ is effective and $h_D(N_G(H)) = 0$, there is an associated real $N_G(H)$-representation $V_D$ associated to $D$. Then by the Mackey formula, $V' := \Res^G_{N_G(H)}\Ind^G_{N_G(H)}(V_D)$ admits $V_D$ as a direct summand.
        Let $V' := V_D \oplus U$, and set $D' = D \otimes C_U\cong C_{V'}$. By $p$-biset functoriality of the dimension homomorphism (see e.g. \cite{BoYa07}), $D' \cong \Res^G_{N_G(H)}(C_{\Ind^G_{N_G(H)} V_D})$. Set $f' = f \otimes \iota^H_{C_U}$, then we have $\check\Psi^H_{N_G(H)}(f') = \check\Psi^H_{N_G(H)}(f)$ up to scalar by properties of forerunners, after using the identification $\check\Psi^H_{N_G(H)}(C_U) \cong k[h_{C_U}(H)]$.

        Let $C^?:= \tilde{C}((S^{\R G\ominus \R})^?)$ denote the effective endotrivial associated to the regular representation (shifted such that $h_C(G) = 0$) and set $C := C^1$. From \Cref{prop:rep_sphere_props_transferred}, we have that
        \[\Psi^H(C_{h_C(H)})\cong C_{h_C(H)}^H \cong k[N_G(H)/H]^{\oplus i}\] for some $i \in \N$. In particular, for any permutation basis of $\Psi^H(C_{h_C(H)})$, every basis element has stabilizer $H$. Because $C^?$ is a chain complex of $k\Gamma_G$-modules, we have an inclusion $C^H \to C^1$ of $kN_G(H)$-modules. Consequently, we have an inclusion $\tilde\iota^H_C \colon k \to \Res^G_{N_G(H)} C_{h_C(H)}$ of $kN_G(H)$-modules as in \Cref{thm:construction}. We regard $\tilde\iota^H_C\colon k \to \Res^G_{N_G(H)}(C)[-h_C(H)]$ as a morphism of $kN_G(H)$-chain complexes, which satisfies that $\Psi^H_{N_G(H)}(\tilde\iota^H_C)$ is a quasi-isomorphism. Moreover, the target of $\tilde\iota^H_C$ is a direct sum of permutation modules isomorphic to $k[N_G(H)/H]$.

        We define \[\hat{f} := \tr^G_{N_G(H)}(f' \otimes \tilde\iota^H_C) = \sum_{g \in [G/N_G(H)]} g(f' \otimes \tilde\iota^H_C),\] see e.g. \cite[Definition 3.6,2]{Be981}. Note the transfer is well-defined since both $f$ and $\tilde\iota^H_C$ are maps between chain complexes restricted to $N_G(H)$. We claim
        \[\check\Psi^H_G(\hat{f}) = \check\Psi^H_{N_G(H)}(f') = \check\Psi^H_{N_G(H)}(f)\in \Hom_{\on{D}_b(k[N_G(H)/H])}(k, k[h_D(H)+s])\]
        up to a possible scalar, with the latter equality already shown. Indeed, we first observe that
        \[\Psi^H_G(\hat{f}) \cong \Psi^H_{N_G(H)}(f') \otimes \Psi^H_{N_G(H)}(\tilde\iota^H_C).\]
        This follows since $\Psi_{N_G(H)}^H(\tilde\iota^H_C)$ has image in degree $0$ contained inside a direct sum of permutation modules of the form $k[N_G(H)/H]$. Therefore, in the trace sum, $g(f'\otimes \tilde\iota^H_C) $ for every $g \not\in N_G(H)$ has image containing no $H$-fixed cells. Now, $\check\Psi^H_{N_G(H)}(\tilde\iota^H_C)$ is a quasi-isomorphism, and up to scalar identifies with $\id_k$ in $\on{D}_b(k[N_G(H)/H])$. The claim follows, and we may scale $\hat{f}$ to identify $\check\Psi^H_G(\hat{f})$ with $\check\Psi^H_{N_G(H)}(f')$ on the nose, rather than up to scalar. Since $\hat{f}$ is a homomorphism from $k$ to a shift of $C_{\Ind^G_H V} \otimes C$, it is a homogeneous element of $\tH(G)$, and we are done.

    \end{proof}

    We are ready to prove injectivity.

    \begin{theorem}\label{thm:injectivity}
        Let $G$ be a $p$-group. The comparison map \[\on{comp}_G\colon  \Spc(\catK(G)) \to \Spec^h(\tH(G))\] of \Cref{prop:comparisonmap} is injective.
    \end{theorem}
    \begin{proof}
        Let $\catP = \catP(H_1,\mathfrak{p})$ and $\catQ = \catP(H_2, \mathfrak{q})$ in $\Spc(\catK(G))$ satisfy $\on{comp}_G(\catP) = \on{comp}_G(\catQ)$ in $\Spec^h(\tH(G))$. This implies that $\catP \in \open(f)$ if and only if $\catQ \in \open(f)$ for every $f \in \tH(G)$. In particular, for every effective endotrivial $C \in \calB(G)$ and subgroup $H \leq G$, we have that $\catP \in \open(\iota^H_C)$ if and only if $\catQ \in \open(\iota^H_C)$.

        Suppose for contradiction that $H_2 \not\leq_G H_1$, and let $C := C_{\R G}$ denote the effective endotrivial associated to the regular representation. In this case we have $h_C(H) = [G:H]$ for any $H \leq G$, so if $|H_1| \neq |H_2|$, then $\iota_C^{H_1}$ and $\iota_C^{H_2}$ have images (as morphisms of chain complexes) in different homological degree, and since $H_2 \not\leq_G H_1$, $\Psi^{H_2}(\iota_C^{H_1})$ is zero. Otherwise if $|H_1| = |H_2|$ but $H_1 \neq_G H_2$, \Cref{lem:samesizesubs} implies both $\Psi^{H_2}(\iota_C^{H_1})$ is zero and $\Psi^{H_1}(\iota_C^{H_2})$ is zero.

        Now, consider the map $\check{\psi}^{H_2}\colon  V_{\Weyl{G}{H_2}} \hookrightarrow \Spc(\catK(G))$. We have:

        \begin{align*}
            (\check{\psi}^{H_2})\inv(\open(\iota^{H_1}_C))&= \open(\cone(\check{\Psi}^{H_2}(\iota^{H_1}_C)))\\
            &= \open(0\colon  k  \to \check{\Psi}^{H_2}(C))\\
            &= \emptyset
        \end{align*}

        Therefore, $V_{\Weyl{G}{H_2}}\cap \open(\iota^{H_1}_C) = \emptyset$ in $\Spc(\catK(G))$. On the other hand, $V_{\Weyl{G}{H_1}} \subseteq \open(\iota^{H_1}_C)$, since $\iota^{H_1}_C$ is invertible in $\catK(G)/\catM(H_1)$. Therefore, $\catP \in\open(\iota^{H_1}_C) $ but $\catQ \not\in \open(\iota^{H_1}_C)$, a contradiction. Thus $H_2\leq_G H_1$, and by symmetry $H_1 =_G H_2$. Set $H := H_1.$

        If $H \trianglelefteq G$, then the composition \[V_{\Weyl{G}{H}} \xrightarrow{\check{\psi}^{H}} \Spc(\catK(G)) \xrightarrow{\on{comp}_G} \Spec^h(\tH(G))\] is injective from \Cref{prop:diagonalinjective}. Otherwise, we have the following commutative diagram.

        \begin{figure}[H]
            \centering
            \begin{tikzcd}
                V_{N_G(H)/H} \ar[r, hookrightarrow] \ar[dr, hookrightarrow] & \Spc(\catK(N_G(H))) \ar[r, "\rho^G_{N_G(H)}"] \ar[d, "\on{comp}_{N_G(H)}"]& \Spc(\catK(G)) \ar[d, "\on{comp}_G"]\\
                & \Spec^h(\tH(N_G(H))) \ar[r, "\rho^G_{N_G(H)}"'] & \Spec^h(\tH(G))
            \end{tikzcd}
        \end{figure}

        Here, the diagonal arrow is injective since $H \trianglelefteq N_G(H)$, the top row is injective since $\Spc(\catK(G)) = \bigsqcup_{H \in \Sub_p(G)/G} V_{\Weyl{G}{H}}$, and the composition of the above maps is precisely $\check{\psi}^H$. Therefore, to show injectivity, it suffices to show the bottom composition of maps \[V_{N_G(H)/H} \hookrightarrow \Spec^h(\tH(N_G(H))) \to \Spec^h(\tH(G))\] is injective.
        Let $\frakp, \frakq \in \Spec^h(\on{H}^\sbull(N_G(H)/H)) \cong V_{N_G(H)/H}$ be distinct, and let $\frakp', \frakq' \in \Spec^h(\tH(N_G(H)))$ denote their (distinct) images in $\Spec^h(\tH(N_G(H)))$. Then without loss of generality, there exists a homogeneous $f \in \tH(N_G(H))$ such that $\check\Psi^H(f) \in \frakp$ but $\check\Psi^H(f)\not\in\frakq$ by \Cref{rem:equivimage}. Then by \Cref{prop:newhardlemma}, there exists a homogeneous $\hat{f}\in\tH(G)$ such that $\check\Psi^H(\hat{f}) \in \frakp$ but $\check\Psi^H(\hat{f}) \not\in\frakq$, and consequently, applying \Cref{rem:equivimage} again gives that $\rho^G_{N_G(H)}(\frakp') \neq \rho^G_{N_G(H)}(\frakq')$, as desired.

    \end{proof}

    \begin{remark}
        Note that the restriction \[\rho_{N_G(H)}^G\colon \Spec^h(\tH(N_G(H))) \to \Spec^h(\tH(G))\] need not be injective in general, as fusion phenomena will occur, c.f. \cite[Corollary 4.9]{BG25}.

        We remark that for the proofs of  \Cref{prop:newhardlemma} and \Cref{thm:injectivity}, the primary endotrivial required was the endotrivial associated to the regular representation, $C_{\R G}$. This suggests that a `reduced' form of twisted cohomology, in which we only twist by this endotrivial, may possibly be sufficient to obtain these same results. We will analyze this construction in future work.
        Unfortunately, explicitly writing down this endotrivial remains challenging in practice.
    \end{remark}

    As an immediate corollary, we strengthen \Cref{prop:diagonalinjective} and \Cref{prop:dbcommdiagram}.

    \begin{corollary}
        Let $H\leq G$ be any subgroup. Then the commutative square

        \begin{figure}[H]
            \centering
            \begin{tikzcd}
                V_{\Weyl{G}{H}} = \Spc(\on{D}_b(k[\Weyl{G}{H}])) \ar[r, "\check{\psi^H}"] \ar[d, "\on{comp}_{\on{D}_b(k[\Weyl{G}{H}])}"'] & \Spc(\catK(G)) \ar[d, "\on{comp}_G"]\\
                \Spec^h(\on{H}^\sbull(\Weyl{G}{H})) \ar[r, hookrightarrow] & \Spec^h(\tH(G))
            \end{tikzcd}
        \end{figure}
        satisfies that the diagonal is injective. Moreover, \[\check\psi^h\colon  \Spec^h(\on{dH}^\sbull(G)) \to \Spec^h(\tH(G))\] is injective with dense image.
    \end{corollary}

    \begin{remark}
        The comparison map is an injective map of spectral spaces, therefore Stone duality implies the pullback \[\Omega(\Spec^h(\tH(G))) \to \Omega(\Spc(\catK(G))),\, U \mapsto \comp_G\inv(U)\] is an epimorphism in $\textsf{CFrm}$, where $\Omega(X)$ denotes the coherent frame of open subsets of a spectral topological space $X$. However, epimorphisms in the category $\textsf{CFrm}$ of bounded distributive lattices need not be surjective. This is in contrast to the category of posets, where epimorphisms are precisely the surjective maps.

        In fact, if the induced map on opens is surjective, then by definition, the submonoid of $\tH(G)$ of effective endotrivials is by definition an ample family \cite[Definition 5.4]{Gal25}, and by \cite[Proposition 5.8]{Gal25}, the comparison map would be an open immersion. Therefore, in order to completely deduce the topology of $\Spc(\catK(G))$, more work is required.

    \end{remark}

    If we know twisted cohomology is Noetherian, we further can show that the comparison map restricts to a homeomorphism upon localization to the opens $U(H)$, and therefore is an open immersion, extending Balmer--Gallauer's \cite[Theorem 15.3]{BG25}. The proof of this is quite short, thanks to machinery recently developed by Sanders \cite{San25} regarding fully faithful tt-functors.

    \begin{corollary}\label{cor:homeo}
        If $\tH(G)$ is Noetherian, then for any $H \leq G$, the comparison map restricts to a homeomorphism \[\comp_G\colon U(H) \cong \Spec^h(\calO^\sbull_G(H)).\] In particular, in this case $\comp_G\colon \Spc(\catK(G)) \to \Spec^h(\tH(G))$ is an open immersion, with image the open subspace of $\Spec^h(\tH(G))$ with closed points \[\comp_G(\catM(H)) = \{f \in \tH(G)\mid \Psi^H(f) \text{ is not a quasi-isomorphism}\}.\]
    \end{corollary}
    \begin{proof}
        \Cref{rmk:commsquares} and \Cref{thm:injectivity} assert an injective map \[\comp_{\catL_G(H)}\colon \Spc(\catL_G(H)) = U(H) \to \Spec^h(\calO^\sbull_G(H))).\]
        Noetherianity of $\tH(G)$ implies Noetherianity of $\calO^\sbull_G(H)$. Indeed, $\tH(G)[S_H\inv]$ is Noetherian, and one has an order-preserving embedding from the ideals of $\calO^\sbull_G(H)$ into the set of homogeneous ideals of $\tH(G)[S_H\inv]$ by taking extension of the inclusion $\calO^\sbull_(G)(H) \to \tH(G)[S_H\inv]$.

        Now, a theorem of Sanders \cite[Theorem 10.8]{San25} directly implies $\comp_{\catL_G(H)}$ is a homeomorphism, as desired. Thus, globally $\comp_G$ is a homeomorphism onto its image in $\Spec^h(\tH(G))$. It is easy to check the closed points $\catM(H)$ are mapped to the corresponding maximal primes in the image as described in the theorem statement.

    \end{proof}

    The following corollary is immediate.

    \begin{corollary}\label{cor:dirac}
        Let $\calO_G^\sbull$ denote the sheaf of $\Z$-graded rings on $\Spc(\catK(G))$ obtained by sheafifying $U \mapsto \End^\sbull_{\catK(G)|_U}(k)$. If $\tH(G)$ is Noetherian, $(\Spc(\catK(G), \calO_G^\sbull)$ is a Dirac scheme in the sense of \cite{HP23}.
    \end{corollary}

    Of course, the pertinent question now is to deduce when \Cref{cor:homeo} and \Cref{cor:dirac} hold. We do not expect that Noetherianity is necessary to prove \Cref{cor:homeo}, but conjecture it should hold regardless.

    \section{Examples of Noetherianity of twisted cohomology}

    To conclude, we show the twisted cohomology ring is Noetherian in some cases. We continue to assume $G$ is a $p$-group.

    \begin{definition}\label{def:resonto}
        We introduce the following terminology. We refer to homogeneous elements of $\tH(G)$ as \emph{global sections}, following \cite{Gal25}. A global section of the form $f \colon k \to C[s]$ for $C \in \calB(G)$ is \emph{minimally graded}. An effective endotrivial complex $C$ is \emph{res-onto} for $H \leq G$ if it satisfies the following property: for all $s \in \Z$, the homomorphism induced by restriction on the global sections \[\Res^G_H\colon \Hom_{\catK(G)}(k, C[s]) \to \Hom_{\catK(H)}(k, \Res^G_H C[s])\] is surjective. We say $C$ is \emph{res-onto} if it is res-onto for all subgroups of $G$, and a $p$-group $G$ is \emph{res-onto} if all $C \in \calB(G)$ are res-onto.
    \end{definition}

    \begin{example}\label{ex:resontos}
        Being res-onto is restrictive.
        \begin{enumerate}
            \item Balmer--Gallauer compute that all elementary abelian $p$-groups are res-onto; for $u_N$ (c.f. \Cref{ex:exampleendotriv}), \cite[Proposition 12.10]{BG25} proves so. More generally, for any group $G$, the inflated endotrivial $u_N$ is res-onto by the same computation.

            \item The group $C_{p^2}$ is not res-onto. We give an example: let $G = C_{p^2}$ be generated by $\sigma \in G$, and let $C := kG \to kG \to k$ be the effective endotrivial with $d_1$ the augmentation homomorphism and $d_2\colon kG \to kG$ induced by the assignment $1 \mapsto \sum 1 - \sigma$. Set $H = \langle \sigma^p\rangle$. Then $\Hom_{\catK(H)}(k, \Res^G_H C[-1])$ is generated by the map $1 \mapsto \sum_{i \in \F_p} \sigma^{pi}$, which is homotopic to the map $1 \mapsto \sigma^j \sum_{i \in \F_p} \sigma^{pi}$ for any $j \in \F_p$. However, none of these maps belongs to $\im(\Res^G_H)$, but their sum does, which is null homotopic. The issue is that no such map is $G$-stable on the nose, but all are $G$-stable modulo homotopy.

            However, we note that $C_4$ being non-res-onto does not prevent Noetherianity. Work of Georgakopoulos gives that $\tH(C_4; \Z)$ in integral coefficients is finitely generated, in fact in minimally graded global sections \cite{Geo19}.

            \item The dihedral group $D_8$ has a unique faithful irreducible endotrivial (c.f. \Cref{ex:exampleendotriv}) which is not res-onto. In fact, res-onto-ness fails in two distinct ways. We use the presentation \[D_8 = \langle a,b \mid a^4 = b^2 = 1, ab = ba\inv \rangle.\]  First, set $H := \langle a^2, b\rangle$, $K_0 := \langle b\rangle$, $K_1 = \langle ab\rangle$, and $K_2 = \langle a^2b\rangle$; then $K_0, K_2 \leq H$ and are $D_8$-conjugate but $K_1$ is nonconjugate to $K_0, K_2$. Then we have the faithful irreducible \[C = kD_8 \to k[D_8/K_0]\oplus k[D_8/K_1] \to k\] and it follows either by direct computation or by checking h-marks that \[\Res^G_H C \cong kH \to k[H/K_0] \oplus k[H/K_2] \to k.\] Moreover, one may compute that \[\Res^G_H \iota_C^{K_0} = \iota_{\Res^G_H C}^{K_0} + \iota_{\Res^G_H C}^{K_2} \text{  and  } \Res^G_H \iota_C^{K_1} = 0.\] In particular, $\iota_{\Res^G_H C}^{K_0}$ and $\iota_{\Res^G_H}^{K_2}$ are not in the image of restriction. This can be viewed as a fusion phenomenon: $K_0$ and $K_2$ are $D_8$-conjugate but not $H$-conjugate, so $\iota_{\Res^G_H C}^{K_0}$ and $\iota_{\Res^G_H C}^{K_2}$ are not $G$-stable, and cannot be in the image of restriction.

            Additionally, let $H' \leq G$ be the unique subgroup isomorphic to $C_4$. Then $\Res^G_{H'} C_\partial$ is isomorphic to the endotrivial given in the previous example when $p = 2$, and consequently, $C_\partial$ is not res-onto for $Z(D_8)$. On the other hand, explicit computation deduces that $C_\partial$ is res-onto for all other subgroups of $D_8$.

        \end{enumerate}

    \end{example}

    \begin{notation}
        Let $C \in \catK(G)$. We set \[\calF_s^G(C) = \calF_s(C) := \Hom_{\catK(G)}(k, C[-s])\] for shorthand to denote the global sections of $C$. We note that $\Gamma^G_s(C)$ would perhaps be more standard, however $\Gamma_G$ already denotes the orbit category of $G$.
    \end{notation}

    \begin{prop}\label{lem:resontoformofmorphisms}
        Let $C$ be an effective endotrivial complex of $kG$-modules and $H\leq G$. The following diagram
        \begin{figure}[H]
            \centering
            \begin{tikzcd}
                & \calF^G_s(C) \ar[dl, "\Res^G_H"] \ar[dr, "\cdot \on{coaug}"']\\
               \calF^H_s(\Res^G_H C) \ar[rr, "\on{ad}"] \ar[ur, shift left=1.7, "\on{tr}^G_H"] & & \calF^G_s( C \otimes k[G/H]) \ar[lu, "\id \otimes \on{aug}"', shift right = 1.7]
            \end{tikzcd}
        \end{figure}
        commutes in the following sense: \[\cdot \on{coaug} = \on{ad} \circ \Res^G_H\] and \[\on{tr}^G_H = (\id \otimes \on{aug})\circ \on{ad}.\] Here, $\on{ad}$ is an isomorphism arising from the induction/restriction adjunction and Frobenius reciprocity. In particular, $C$ is res-onto for $H$ if and only if $\on{coaug}$ is surjective, that is, every element of $\Hom_{\catK(G)}(k, C[s] \otimes k[G/H])$ is homotopic to a simple tensor $f \otimes \on{coaug}$.
    \end{prop}
    \begin{proof}
        The key point is that $\on{ad}$ sends the $H$-linear global section $f_s \colon k\to \Res^G_H C_s$ in homological degree $s$ to the $G$-linear map induced from the assignment \[k \ni a \mapsto \sum_{g\in [G/H]} gf_s(a) \otimes gH \in C_s \otimes k[G/H].\] Its inverse $\on{ad}\inv$ returns the identity coset in the sum. From these, verifying the identities is routine.

    \end{proof}

    The following is a generalization of \cite[Lemma 12.12]{BG25}.
    \begin{prop}\label{prop:res_onto_case}
        Let $D \in \catK(G)$ contain no contractible direct summands and arise from a chain complex of $k\Gamma_G$-modules, and let $C \in \Pic(\catK(G))$ be indecomposable and effective. Write $D_i = \bigoplus_{j \in J}k[G/H_{i,j}]$, then by \Cref{prop:ctrctablesummandoforbitcc} we have a coaugmentation homomorphism $\gamma_{i,j}\in \calF_i^G(D)$ with image the socle of $k[G/H_{i,j}]$.

        Suppose for every subgroup $H_{i,j}$ appearing as an isotropy subgroup in $D$ that the restriction \[\Res^G_{H_{i,j}}\colon\calF^G_s(C) \to \calF^{H_{i,j}}_s(\Res^G_{H_{i,j}} C)\] is surjective for every $s$. Then \[\calF^G_s(C \otimes D) = \sum_{i,j} \calF^G_{s-i}(C)\gamma_{i,j}.\]
    \end{prop}
    \begin{proof}
        We have for every $i$ a triangle arising from brutal truncation: \[D_i[i] \to D_{\geq i} \to D_{\geq i+1} \to D_i[i+1],\] noting each truncation again contains no contractible direct summands, and consequently an exact sequence
        \[\calF_{s-i}^G(C \otimes D_i)\to \calF^G_s(C \otimes D_{\geq i}) \to \calF^G_s(C \otimes D_{\geq i+1}) \to \calF^G_{s - (i+1)}(C \otimes D_i).\]
        For every transitive summand in $D$ the induction-restriction adjunction gives \[\calF^G_s(C \otimes k[G/H]) \cong \calF^H_s(\Res^G_H C), \] and under this isomorphism the class $f \otimes \on{coaug}$ corresponds to $\Res^G_H f$ by \Cref{lem:resontoformofmorphisms}. Surjectivity implies every class on the left is of this form.

        Now, we induct on the length of $D$, following the strategy of \cite[Lemma 12.12]{BG25}. The base case of $D$ being concentrated in one degree follows from above, via the assumption that restriction is surjective for every $H_{i,j}$. Now, given a $f \in \calF^G_s(C \otimes D_{\geq i})$, its image in $f \in \calF^G_s(C \otimes D_{\geq i+1})$ can be written as a product as specified. This image can be lifted to the `same' element in $\calF^G_s(C \otimes D_{\geq i})$ because each global section $\gamma$ is defined for both truncations $D_{\geq i+1}$ and $D_{\geq i}$. Subtracting $f$ with this term gives an element in $\calF^G_{s - i}(C \otimes D_i)$, and taking its image in $\calF^G_s(C \otimes D_{\geq i})$ by multiplying with each $\gamma_{i,j}$ using \Cref{lem:resontoformofmorphisms} expresses $f$ via the global sections $\gamma$ and global sections of $C$, as desired.

    \end{proof}

    \subsection{Noetherianity for $D_8$} Res-onto-ness is not a necessary condition for Noetherianity. We prove $\tH(G)$ is Noetherian for $D_8$, the dihedral group of 8 elements.

    \begin{lemma}\label{lem:gstables}
        Set $G = D_8$ and $H \leq G$ isomorphic to $V_4$. Let $C \in \Pic(\catK(G))$ be an effective endotrivial and $s \in \Z$. Then
        \begin{align*}
            \calF_s^H(\Res^G_H C)^G &= \im\big(\Res^G_H\colon \calF_s^G(C) \to\calF_s^H(\Res^G_H C)\big) \\
            &= \ker\big(\tr^G_H\colon \calF_s^H(\Res^G_H C)\to  \calF_s^G(C)\big)
        \end{align*}
    \end{lemma}
    \begin{proof}
        We have a set of inclusions \[\im(\Res^G_H)\subseteq \ker(\tr^G_H) \subseteq \calF_s^H(\Res^G_H C)^G,\] so it suffices to show every $G$-stable element belongs to $\im(\Res^G_H)$. This computation is performed explicitly on twisted cohomology, where we can describe the graded (homogeneous) subring $\tH_h(H)^G$ explicitly. We have $\calB(H) = \{D_1,D_2,D_3\}$, with $D_i = k[H/K_i] \to k$, with $K_i \leq H$ a subgroup of order 2. If we set $K_3 = Z(G)$, then $G$ acts by swapping $D_1, D_2$. From \cite{HS25}, we have, using the notation of \cite{BG25}, a minimal presentation from \Cref{rmk:bgmorphisms}
        \[\tH(H) = k[ a_{K_1}, b_{K_1}, a_{K_2}, b_{K_2}, a_{K_3}, b_{K_3}]/\langle a_{K_1}b_{K_2}b_{K_3} + b_{K_1}a_{K_2}b_{K_3} + b_{K_1}b_{K_2}a_{K_3}\rangle.\] We show the graded subring of \emph{homogeneous} $G$-stable elements is generated by the following five elements:
        \[\tH_h(H)^G = \langle a_{K_3}, b_{K_3}, a_{K_1}a_{K_2}, b_{K_1}b_{K_2}, a_{K_1}b_{K_2} + b_{K_1}a_{K_2}\rangle.\] Indeed, $G$ exchanges $a_{K_1}$ with $a_{K_2}$, $b_{K_1}$ with $b_{K_2}$.
        Therefore, $\tH_h(H)^G$ is contained in the graded subring $\tH(H)'$ of $\tH(H)$ which coincides with $\tH(H)$ degrees parameterized by the endotrivials written as products of $D_1D_2$ and $D_3$, and is 0 elsewhere. Using \Cref{prop:res_onto_case} (as $V_4$ is res-onto) gives that \[\tH(H)' =\langle a_{K_3}, b_{K_3}, a_{K_1}a_{K_2}, b_{K_1}b_{K_2}, a_{K_1}b_{K_2}, b_{K_1}a_{K_2}\rangle.\] However, the $G$-action exchanges $ a_{K_1}b_{K_2}$ and $b_{K_1}a_{K_2}$.  Every homogeneous element of $\tH(H)'$ can be written (not necessarily uniquely) in the form $f = f' + f''\cdot a_{K_1}b_{K_2}$ where $f', f''$ are homogeneous and belong to the subring \[\tH(H)'':= \langle a_{K_3}, b_{K_3}, a_{K_1}a_{K_2}, b_{K_1}b_{K_2}, a_{K_1}b_{K_2} + b_{K_1}a_{K_2}\rangle.\] Since $f', f''$ are $G$-stable, it follows that $f$ is $G$-stable if and only if $f''\cdot (a_{K_1}b_{K_2} + \cdot b_{K_1}a_{K_2}) = 0$. As $\tH(H)$ is an integral domain (as the one defining relation is an irreducible polynomial), we conclude $f'' = 0$. Thus $\tH_h(H)^G = \tH(H)''$, as desired.
        Note that we have non-homogeneous $G$-stable elements $a_{K_1} + a_{K_2}$ and $b_{K_1} + b_{K_2}$ as well.

        Now, each of these generators is in the image of $\Res^G_H\colon \tH(G) \to \tH(H)^G$, with the latter three terms coming from forerunners for $C_\partial$, and the first two coming from forerunners for the endotrivial $C_L := k[G/L]\to k$, with $L \leq G$ the cyclic subgroup of order 4. It remains to show that every element of $\calF_s^H(\Res^G_H C)^G$ is mapped onto by an element of $\calF^G_s(C)$ via restriction (as opposed to from some other twist). Set $C_H := k[G/H] \to k$ and set $C_{H'} := k[G/H'] \to k$ where $H'\leq G$ is the other subgroup isomorphic to $V_4$. We have \[\Res^G_H C_H \cong k[1], \quad \Res^G_H C_{H'} \cong \Res^G_H C_L \cong D_3,\quad \Res^G_H C_\partial \cong D_1D_2.\]
        First, if $C$ is one of the complexes $C_H, C_{H'}, C_L,$ or $C_\partial$, then this is easily verified. Then more generally, if $C = C_H^wC_{H'}^xC_L^yC_\partial^z$ and $f \in \calF_s^H(\Res^G_H C )^G$, then $\Res^G_H C = (D_1D_2)^{z}D_3^{x+y}[w]$. Writing $f$ as a polynomial in the generators of $\tH_h(H)^G$, the generators graded in $D_1D_2$ lift to $C_\partial$, and the generators graded in $D_3$ can be given lifts distributed between the $x$ copies of $C_{H'}$ and $y$ copies of $C_L$. Finally, $C_H$ accounts for the shift $[w]$, and we conclude $f$ can be lifted to $\calF_s^G(C)$.

    \end{proof}

    \begin{remark}\label{rmk:mistake}
        In a previous version of this paper, \Cref{lem:gstables} was incorrectly asserted to hold for $G$ a $p$-group and $H \leq G$ a subgroup satisfying $[G:H] = p$. The analogue of this statement is easily shown to hold over the category of chain complexes, but up to homotopy, the claim fails. It was also previously incorrectly asserted that all Dedekind $p$-groups are res-onto. For instance, \Cref{ex:resontos}(b) is a counterexample to both of these claims.

    \end{remark}

    \begin{theorem}\label{thm:d8noeth}
        Set $G = D_8$. Then $\tH(G)$ is generated as graded $k$-algebra by its minimal morphisms. In particular, $\tH(G)$ is Noetherian.
    \end{theorem}
    \begin{proof}
        Explicit computation shows that the only obstruction to $G$ being res-onto comes from the faithful irreducible endotrivial $C_\partial$ (see \Cref{ex:resontos}), when restricted to the two subgroups of $D_8$ isomorphic to $V_4$ and $Z(D_8)$, however $Z(D_8)$ does not appear as an isotropy subgroup in any irreducible effective endotrivial of $D_8$, so this case can be ignored. Given $C = C_1 \otimes \cdots \otimes C_n$ a product of elements in $\calB(G)$, we proceed by induction on $n$.
        If either no $C_i$ is isomorphic to $C_\partial$, or some $C_i$ is isomorphic to $C_\partial$ but no $C_i$ is isomorphic to the endotrivial $k[G/H]\to k$ with $H \cong V_4$, then by induction and \Cref{prop:res_onto_case}, every $f \in \calF_\bullet^G(C)$ can be written as a polynomial of global sections. Indeed, in either case we can choose in the lemma for $D$ the direct summand of $C_1 \otimes \cdots \otimes C_{n-1}$ with no contractible summands and $C = C_n$. Note that in the latter case, $Z(D_8)$ never appears as an isotropy subgroup.

        Otherwise, we may assume without loss of generality that $C_n = k[G/H] \to k$ with $H \cong V_4$. Let $s \in \Z$, set $D := C_1 \otimes \cdots \otimes C_{n-1}$ and $C:= k[G/H] \to k$ for $H \cong V_4$. We again have an exact triangle
        \begin{figure}[H]
            \centering
            \begin{tikzcd}
                k[0] = &  \cdots 0 \ar[r]  & 0  \ar[r] \ar[d] & k \ar[d, equals, "\iota_C^G"] \\
                C = &  \cdots 0 \ar[r ]  & k[G/H] \ar[r, "\on{aug}"] \ar[d, equals] & k \ar[d] \\
                k[G/H][1] := & \cdots 0 \ar[r]  & k[G/H] \ar[r] \ar[d, "\on{aug}"] & 0 \ar[d] \\
                k[1] = & \cdots 0 \ar[r]& k \ar[r] & 0
            \end{tikzcd}
        \end{figure}
        Tensoring with $D$ and applying $\Hom_{\catK(G)}(k, -[-s])$ gives the exact sequence
        \[\calF^G_s(D) \xrightarrow{\cdot\iota_C^G} \calF^G_s(D \otimes C)\xrightarrow{\,} \calF^G_{s-1}( D \otimes k[G/H])\xrightarrow{\id \otimes \on{aug}}  \calF^G_{s-1}(D),\]

        and the final three terms identify under adjunction and \Cref{lem:resontoformofmorphisms} as
        \[\calF_s^G(D \otimes C) \xrightarrow{\phi}\calF_{s-1}^H(\Res^G_H D) \xrightarrow{\tr^G_H} \calF_{s-1}^G(D).\]
        Let $f \in \calF_s^G(D \otimes C)$ and let $f'$ denote its image in $\calF_{s-1}^H(\Res^G_H D)$. By induction, $f'$ can be expressed as a polynomial in global sections belonging to $\calF_\bullet^H(\Res^G_H C_i)$. However, because $D_8$ is not res-onto for subgroups isomorphic to $V_4$, we cannot automatically lift the same polynomial to one in $\calF_s^G(D\otimes C)$ as in \Cref{prop:res_onto_case}.

        Since the sequence is exact, $\im(\phi) = \ker(\id\otimes \on{aug})$. \Cref{lem:gstables} and \Cref{lem:resontoformofmorphisms} demonstrate that in fact, $\ker(\on{tr}^G_H) = \im(\Res^G_H)$ with $\im(\Res^G_H)$ equivalently the $G/H$-stable global sections. Equivalently by \Cref{lem:resontoformofmorphisms}, $\ker(\id\otimes \on{aug}) = \im(\cdot \on{coaug})$. In particular, the image of $\phi$ consists precisely of the $G$-stable morphisms, i.e., the (homogeneous part of the) ring of invariants.
        By induction, $f'$ is a $G$-stable polynomial in minimally graded elements.
        Then \Cref{lem:gstables} asserts that $f'\in \calF^H_{s-1}(\Res^G_H D)$ lifts to a global section $f'' \in \calF^G_{s-1}(D)$ which satisfies $f - \iota^1_C f'' \in \ker(\phi)$, hence $f- \iota^1_Cf''\in \im(\cdot\iota^G_C)$. Then $f - \iota^1_C f'' = \iota^G_C f'''$ for some $f''' \in \calF_s^G(D)$, and by induction, $f'''$ is also polynomial in minimally graded global sections. Thus $f = \iota^1_C f'' + \iota^G_C f'''$, and we conclude $\tH(G)$ is generated by the minimally graded global sections, as desired.

    \end{proof}

    \begin{remark}
        Determining relations between the generators of $\tH(D_8)$ remains open; while one can determine relations for the global sections of the length two endotrivials arising from inflation from $C_2$ by applying $\Psi^{Z(D_8)}$, the relations between those global sections and those of $C_\partial$ remain unclear.

        We expect that a similar proof strategy should work for extraspecial $p$-groups for $p$ odd, but the computation is much more intensive. These extraspecial $p$-groups have a unique faithful irreducible endotrivial of length $2p+1$, and one sees that the obstruction to being res-onto for these groups is analogous to the case of $D_8$. One key issue is that since $[G:H] > 2$ in odd characteristic, one no longer immediately obtains the equality $\ker(\on{tr}^G_H) = \im(\Res^G_H)$ as before, which is crucial for the strategy we use. Although explicit computation suggests that this equality should hold, we do not have a complete description of the relations of $\tH(G)$ for elementary abelian $p$-groups $G$ with $p$ odd (as \cite{HS25} considers only the case of $p = 2$).
    \end{remark}

    \subsection{Noetherianity for Dedekind $p$-groups}

    We turn to the case of Dedekind $p$-groups. It is well-known that either a Dedekind $p$-group is abelian or isomorphic to $Q_8 \times E$ for an elementary abelian 2-group $E$. Consequently, subgroups of Dedekind $p$-groups are again Dedekind.
    The following proposition can be deduced from the representation theory of Dedekind $p$-groups, but we give a quick proof.

    \begin{prop}\label{prop:dedekind_endotrivs}
        Let $G$ be a Dedekind $p$-group. Then every indecomposable effective endotrivial $C\in \calB(G)$ is of the form \[\Inf^G_{G/N}D,\] where $D$ is a truncated projective resolution of the trivial $k[G/N]$-module. In particular, $G/N$ is cyclic or isomorphic to $Q_8$.
    \end{prop}
    \begin{proof}
        Let $V$ denote the irreducible real representation associated to $C$. If $v \in V$ is not fixed by $g \in G$, then since $\langle g\rangle \trianglelefteq G$, $V^{\langle g\rangle}$ is again a real representation of $G$. However, since  $V^{\langle g\rangle} \subsetneq V$, by irreducibility $V^{\langle g\rangle} = 0$. Therefore, setting $N = \ker(V)$, the action of $G/N$ (which remains Dedekind) on $S^V$ has two fixed points associated to the poles, and consequently, $C_{V}$ consists entirely of free $k[G/N]$-modules in all nonzero degrees. The result follows. The final statement follows from the well-known fact that the only $p$-groups with periodic cohomology have $p$-rank one.

    \end{proof}

    \begin{notation}
        Given an indecomposable endotrivial $C$ and $i \in \N$, we let $C(i)$ denote the unique up to isomorphism indecomposable endotrivial direct summand of $C^{\otimes i}$ (with $C(0) = k$).
    \end{notation}

    \begin{remark}\label{rem:truncated_per_resos}
        If $C \in \calB(G)$ is indecomposable, it is a minimally truncated effective endotrivial of length $d+1$, that is, its maximal nonzero homological degree is degree $d$. In this case, if $C_d = k[G/N]$, then $G/N$ is a periodic Dedekind $p$-group of period $d$. One may compute via h-marks that $C(i)$ is also a truncated effective endotrivial given by a truncated periodic resolution of length $di + 1$, i.e., its maximal homological degree is degree $di$. Therefore, the global sections for these endotrivials also satisfy periodicity conditions, which can already be observed in the cohomology ring $\on{H}^\sbull(G/N;k) \cong \calO_{G/N}^\sbull(N/N)$.

    \end{remark}

    \begin{prop}\label{prop:dedekind_actions_of_iotas}
        Let $G$ be a Dedekind $p$-group and $C \in \calB(G)$ have length $d+1$, i.e., $h_C(1) = d$. Then the bigraded $k$-subalgebra \[\tH_C(G):= \bigoplus_{i \in \N, s\in \Z}\calF^G_s(C(i))\subseteq \tH(G)\] is generated by the morphisms in degrees $\Hom_{\catK(G)}(k, C[s])$ for $s \in \Z$. In particular, $\tH_C(G)$ is finitely generated. Moreover, if $0 \leq s \leq di$, the following hold:
        \begin{enumerate}
            \item We have an isomorphism $\calF^G_s(C(i)) \to \calF^G_{s+d}(C(i+1))$ induced by $\cdot \iota^1_C$.
            \item We have an isomorphism $\calF^G_s(C(i)) \to \calF^G_s(C(i+1))$ induced by $\cdot \iota^G_C$.
        \end{enumerate}
    \end{prop}
    \begin{proof}
        Finite generations follows from \Cref{prop:res_onto_case} since every complex only contains isotropy subgroups $N$ and $G$ (after inflating) with $G/N$ cyclic or $Q_8$, and one easily checks that both these groups satisfy res-onto-ness for $C$.

        For the latter, note that it suffices to assume $G$ is a periodic group of period $d$ and $C$ is a truncated periodic resolution from \Cref{prop:dedekind_endotrivs} and \Cref{rem:truncated_per_resos}. These also imply we have isomorphisms of $k$-vector spaces $\calF^G_s(C(i)) \cong \calF^G_{s+d}(C(i+1))$ and $\calF^G_s(C(i)) \cong \calF^G_s(C(i+1))$. It remains to show that multiplication by the forerunners induces these isomorphisms. This can already be checked on group cohomology (i.e., in $\calO^\sbull_G(1)$), as if $\calF^G_s(C(i))$ is nonzero, localization induces an isomorphism \[\calF^G_s(C(i)) \cong \on{H}^{di - s}(G;k).\] In this case, multiplication by $\iota^1_C$ is the identity map (as $\iota^1_C = 1$ in $\calO^\sbull_G(1)$), and multiplication by $\iota^G_C$ is the periodicity isomorphism \[\on{H}^{di - s}(G;k) \cong \on{H}^{d(i+1)-s}(G;k).\]

    \end{proof}

    The following alternative characterization of $\tH(G)$ will be crucial for proving finite generation.

    \begin{prop}\label{lem:alt_description_of_th}
        Let $G$ be a $p$-group. For any endotrivial $C$, let $C(i)$ denote the indecomposable direct summand of $C^{\otimes i}$ (with $C(0) = k$). We define an infinite-dimensional dg-module $\calM$ as follows: \[\calM := \bigoplus_{q \in \N^{\calB(G)}} \bigotimes_{C \in \calB(G)} C(q(C)) = \bigotimes_{C \in \calB(G)}\bigoplus_{i \in \N} C(i).\] Then as multigraded $k$-modules, $\on{H}_\sbull(\calM^G) \cong \tH(G)$, where $\calM^G$ denotes the $G$-invariants of $\calM$.
    \end{prop}
    \begin{proof}
        Indeed, for any $C$ a product of elements in $\calB(G)$ and $s \in \Z$, $Z(C_s^G)$ is precisely the group of global sections $k \to C_s$ (not modulo homotopy), and a global section is null homotopic if and only if its image is also in the the image of the differential $d_{s+1}^G$ of $C^G$. Explicitly, the isomorphism sends a global section $f\colon k \to C[-s]$ in $\calF_s^G(C) \subset \tH(G)$ to the element $f(1) \in C_s$, which one may check is a well-defined element of $\on{H}_s(C^G)$.

    \end{proof}

    Note that by \Cref{prop:identifyingkhoms}, if an endotrivial is indecomposable (by taking the unique up to isomorphism direct summand), then in fact there are no coboundaries; that is, for every $C(q(C))$, the differentials of the $G$-invariants are trivial.

    \begin{definition}\label{def:S_const}
        Let $G$ be a $p$-group. For $C \in \calB(G)$, set \[\calM_C := \bigoplus_{i \in \N}C(i).\] Then $\calM= \bigotimes_{C\in \calB(G)} \calM_C$. Then $\on{H}_\sbull(\calM_C^G)$ corresponds via the isomorphism in \Cref{lem:alt_description_of_th} to the bigraded subalgebra \[\tH_C(G):= \bigoplus_{i \in \N, s\in \Z} \calF^G_s(C(i)) \subseteq \tH(G).\]

        If $G$ is Dedekind, then from \Cref{lem:alt_description_of_th} and \Cref{prop:dedekind_actions_of_iotas}, the actions of the forerunners $\iota^G_C, \iota^1_C\in \tH(G)$ induce actions on $\on{H}_\sbull(\calM_C^G)$. We extend this to an action on $\calM_C$ in the following way. Let $\calS_C$ denote the commutative $k$-algebra on two variables, $k[a_C, b_C]$. We define an action on $\calM_C$ as follows. Assume $C$ is a truncated $d$-periodic resolution of the $k[G/N]$-module $k$, then for all $s \leq  di$ we have identifications $C(i)_s = C(i+1)_s$ and for all $0 < t \leq d(i-1)$ we have identifications $C(i)_t = C(i)_{t+d}$. For $m \in C(i)_s$, set $a_C \cdot m := m \in C(i+1)_s$. Then by periodicity, multiplication by $a_C$ is a well-defined chain map on $\calM_C$.

        For $b_C$, we need to be a bit more careful. If $m \in C(i)_s$ with $s \neq 0$, we set $b_C \cdot m := m\in C(i+1)_{s+d}$. In this case, it is clear by periodicity that this is a well-defined chain map. Otherwise, if $s = 0$, then $C(i)_0 = k$ but $C(i+1)_{d}= k[G/N]$ for some $N \leq G$. In this case, for $m \in C(i)_0$, we set $b_C \cdot m := \on{coaug}(m) \in C(i+1)_d$. This again is a well-defined chain map, as the image of $b_C$ sends $m$ to the socle of $C(i+1)_d$, which belongs to $\ker(d_i)$.

    \end{definition}

    \begin{lemma}\label{prop:fin_gen_plus_factor}
        Let $G$ be a Dedekind $p$-group. With the notation as in \Cref{def:S_const}, the following hold for every $C \in \calB(G)$.
        \begin{enumerate}
            \item $\calM_C$ is a finitely generated $\calS_C$-module. Therefore, $\tH_C(G)$ is a finitely generated $\calS_C$-module.
            \item The action of $a_C$ on $C(i)\subset \calM_C$ identifies with $\iota^G_C \otimes \id_{C(i)}$ up to scalar.
            \item The action of $b_C$ on $C(i) \subset \calM_C$ identifies with $\iota^1_C \otimes \id_{C(i)}$ up to scalar.

        \end{enumerate}
    \end{lemma}
    \begin{proof}
        (a) is an easy inductive verification, following from periodicity of $C$. In particular, $\calM_C$ is generated by a $k$-basis of $C(1)$ and a nonzero element of $C(0)$.

        For (b), we have $\Hom_{\catK(G)}(k, C) \cong k$, with generator $\iota^G_C$. By invertibility of $C$, it follows that $\Hom_{\catK(G)}(C(i), C(i+1)) \cong k$ as well, and it is easy to check that $\iota^G_C \otimes \id_{C(i)} \in \Hom_{\catK(G)}(C(i), C(i+1))$ is nonzero. Since the action of $a_C$ is nonzero in $\catK(G)$ as well, (b) follows. (c) follows analogously.

    \end{proof}

    We are now ready to prove finite generation.

    \begin{theorem}\label{thm:resontonoeth}
        Let $G$ be a Dedekind $p$-group. Then $\tH(G)$ is a finitely generated $k$-algebra, in particular Noetherian.
    \end{theorem}
    \begin{proof}

        By \Cref{prop:fin_gen_plus_factor}(a), for every $C \in \calB(G)$, $\calM_C$ is a finitely generated $\calS_C$-module. Therefore by right exactness of the tensor product, $\calM$ is finitely generated over \[\calS:= \bigotimes_{C \in \calB(G)} \calS_C.\] The action of $\calS$ preserves $G$-invariants of $\calM$, so consequently, $\calM^G$ is finitely generated over $\calS$ as well. Again, $\calS$ preserves and reflects cocycles and preserves coboundaries in $\calM^G$, hence $\on{H}_\sbull(\calM^G)$ is finitely generated over $\calS$.

        Therefore, $\tH(G)$ is finitely generated over $\calS$ too. Now, \Cref{prop:fin_gen_plus_factor}(b) implies $\tH(G)$ is finite over the finitely generated subalgebra \[\langle \iota^1_C,\iota^G_C \mid C \in \calB(G) \rangle,\] where the action is given by \Cref{prop:fin_gen_plus_factor}(b,c), whence the result.

    \end{proof}

    \begin{example}
        Dedekind $p$-groups give us examples of twisted cohomology rings which are finitely generated, but not generated by minimally graded morphisms. These arise from the failure of res-onto-ness exhibited in \Cref{ex:resontos}(b). For instance, let $p = 2$ and $G = C_4 \times C_2 = \langle\sigma, \tau\mid \sigma^4 = \tau^2 =  1, \sigma\tau = \tau\sigma\rangle$. Let $H_1 = \langle \tau\rangle$ and $H_2 = \langle \sigma^2\tau\rangle$, then we have two effective endotrivials \[C^1 = k[G/H_1] \to k[G/H_1]\to k\] and \[C^2 = k[G/H_2] \to k[G/H_2] \to k.\] Explicitly, these endotrivials are both inflated from the example of \Cref{ex:resontos}(b). Then $(C^1 \otimes C^2)_2 = k[G/H_1] \oplus k[G/H_2] \oplus k[G/H_1 \times G/H_2]$. Let $C$ denote the unique indecomposable endotrivial direct summand of $C^1 \otimes C^2$. Denote the differentials of $C, C^1, C^2$ by $d_i, d^1_i, d^2_i$ respectively.

        First, we claim $C_2$ contains a direct summand isomorphic to $kG$. Indeed, it is straightforward that $C_0 = k$ and $C_1 = k[G/H_1] \oplus k[G/H_2]$. First, $C$ has a permutation basis satisfying the stabilizers grow property by \Cref{prop:removecontractiblesummands}, and there are no nonzero maps $k[G/H_1] \to k[G/H_2]$ or $k[G/H_2] \to k[G/H_1]$ satisfying stabilizers grow. However, $\ker(d_1)$ contains elements $m$ for which $p_{H_1}(m), p_{H_2}(m) \not\in\ker(d_1)$. Since $C$ is exact in homological degree 1, $m \in \im(d_2)$. Neither $k[G/H_1]$ nor $k[G/H_2]$ can cover $m$, so there must be a projective direct summand of $C_2$ mapping onto it. In particular, $\dim_k \calF^G_2(C) \geq 3$.

        Now, we claim the homomorphism \[\phi\colon \bigoplus_{i + j = 2} \calF_i^G(C^1) \otimes \calF_j^G(C^2) \to \calF_2^G(C^1 \otimes C^2) \cong \calF_2^G(C)\] induced from tensoring cannot be surjective. Let $f \in \calF_1^G(C^1)$ and $g \in \calF_1^G(C^2)$ be nonzero, in which case they span their respective spaces of global sections. One may observe that $f \otimes g$ is null homotopic; this computation follows analogously to the computation in \Cref{ex:resontos}(b) after using adjunction. Therefore, the domain of $\phi$ has $k$-dimension 3 with kernel having $k$-dimension at least 1, but the target has $k$-dimension at least 3. Consequently, the subring of $\tH(G)$ generated by minimally graded morphisms does not contain all of $\calF_2^G(C)$.
    \end{example}

     We emphasize that showing Noetherianity of twisted cohomology is not the only path to showing that the comparison map is an open immersion. However, Noetherianity remains an interesting, and evidently rather difficult, question.

    \bibliography{bib}
    \bibliographystyle{alpha}
\end{document}